\documentclass{amsart}
\usepackage{amssymb}
\usepackage[mathscr]{eucal}
\usepackage{url}
\usepackage[all]{xy} 
\usepackage{graphicx}
\usepackage[new]{old-arrows}
\usepackage[usenames,dvipsnames]{xcolor}
\definecolor{brightcerulean}{rgb}{.26, .41, .88}

\usepackage[ allcolors=brightcerulean, allbordercolors=brightcerulean,
        colorlinks, citecolor=brightcerulean,
        pdfauthor={},
        pdftitle={}
]{hyperref}

\usepackage{tikz, enumerate,caption,stmaryrd,amsfonts,amssymb,systeme,comment,float,mathrsfs, mathtools, bigfoot}
\usetikzlibrary{matrix,positioning,arrows,shapes,chains,calc,automata}
\usetikzlibrary{decorations.pathreplacing}
\usepackage{tikz-cd}
\usepackage{float}

\makeatletter
\let\@wraptoccontribs\wraptoccontribs
\makeatother

\newtheorem{thm}[equation]{Theorem}
\newtheorem{lem}[equation]{Lemma}
\newtheorem{cor}[equation]{Corollary}

\newtheorem{prop}[equation]{Proposition}
\newtheorem{conj}[equation]{Conjecture}

\theoremstyle{definition}
\newtheorem{rem}[equation]{Remark}

\newtheorem{defn}[equation]{Definition}

\newfont{\cyrr}{wncyr10}

\def\Z{\mathbb{Z}}
\def\Q{\mathbb{Q}}
\def\F{\mathbb{F}}

\def\R{\mathbb{R}}
\def\C{\mathbb{C}}

\def\fp{\mathfrak{p}}
\DeclareMathOperator{\HH}{\mathrm{H}}
\newcommand{\et}{\text{\'et}}
\newcommand{\dR}{\mathrm{dR}}

\def\Div{\mathrm{Div}}

\def\Gal{\mathrm{Gal}}
\def\rk{\mathrm{rank}}

\def\coker{\mathrm{coker}}

\def\Symm{\mathrm{Symm}}

\def\Sel{\mathrm{Sel}}

\def\Aut{\mathrm{Aut}}
\def\loc{\mathrm{loc}}

\def\Pic{\mathrm{Pic}}
\def\Nm{\mathrm{Nm}}
\def\Tr{\mathrm{Tr}}
\def\NS{\mathrm{NS}}

\def\onto{\twoheadrightarrow}

\def\vv{\vskip10pt}
\def\ed{\end{document}}

\begin{document}

\title{Ogg's Torsion Conjecture: Fifty years later}
\vskip10pt

\author[J. S. Balakrishnan and B. Mazur (with an appendix by N. Dogra)]{Jennifer S. Balakrishnan and Barry Mazur \\(with an appendix by Netan Dogra)
}

\begin{abstract} Andrew Ogg's mathematical viewpoint has inspired an increasingly  broad array of results and conjectures. His results and conjectures have earmarked fruitful turning points  in our subject, and his influence has been such a gift to all of us.\footnotemark
\vv
 Ogg's celebrated Torsion Conjecture---as it relates to modular curves---can be paraphrased as saying that rational points (on the modular curves that parametrize torsion points on elliptic curves) exist if and only if there is a good geometric reason for them to exist.\footnotemark 
 \vv
 We give a survey of Ogg's Torsion Conjecture and the subsequent developments in our understanding of rational points on modular curves over the last fifty years.
 \end{abstract}
\maketitle
\stepcounter{footnote}
\footnotetext{This paper expands  the 45-minute talk that B.M. gave at the conference at the IAS celebrating the Frank C. and Florence S. Ogg Professorship  in Mathematics. \vv}
\stepcounter{footnote}

\footnotetext{\emph{B.M.: }And here's just one  (tiny) instance of Ogg's jovial and joyful way of thinking:  As Tate and  I recorded in one of our papers \cite{MT}:
\vv
  ``Ogg passed through our town and mentioned that he had discovered a point of order $19$''  on the Jacobian of $X_1(13)$, allowing us to feel that that Jacobian was 
{``not entitled to have''}
 more than $19$ points.\vv}
 \tableofcontents

 \begin{center}
 \begin{figure}[!h]
 \caption{Andrew Ogg, photo by George M. Bergman \\(Oberwolfach Photo Collection)}
  \vskip10pt \includegraphics[width= 0.8\textwidth]{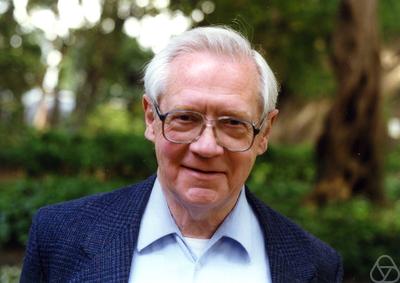}\end{figure}\end{center}
 \part{An overview}\label{overview}

    Torsion in algebraic groups---even if not in that vocabulary---has played a fundamental role since Gauss's  {\it Disquisitiones Arithmeticae} (1801), the structure of roots of unity (torsion in the multiplicative group) being a  central concern in the development of modern number theory.{\footnote{See  Umberto Zannier's expository article  \cite{Zan},
{\it Torsion in algebraic groups and problems which arise.}}}\vv
      Let $K$ be a number field,  and denote by $G_K$ its absolute Galois group, i.e. $$G_K \coloneqq \Gal({\overline{K}}/K).$$ \vv
       A basic question in the arithmetic of abelian varieties over number fields is to classify (up to the natural notion of isomorphism) pairs  $$(A; C \stackrel{\alpha}{\longhookrightarrow}A({\overline{K}}))$$  
       
        \noindent where \begin{itemize}\item $A$ is a (polarized) abelian variety defined over $K$,  \item $C$ is a finite abelian  group with a $G_K$-action, and 
       \item $\alpha$ is a $G_K$-equivariant injection. \end{itemize}
       \vv
      These are the three basic parameters in this general question, and you have your choice of how you want to choose the range of each of them.  For example, you can:
      \vv
       \begin{itemize}\item allow the groups $C$ to run through all cyclic finite groups with arbitrary $G_K$-action; and $A$ to range through all abelian varieties with a specified type of polarization.  Equivalently, you are asking about {\bf $K$-rational cyclic isogenies of abelian varieties}, or \vv \item  restrict to finite groups $C$ with trivial $G_K$-action, in which case you are asking about {\bf $K$-rational torsion points on abelian varieties},\vv \item vary over a class of number fields $K$---e.g., number fields that are of a fixed degree $d$ over a given number field $k$, or
    \vv   \item fix the dimension of the abelian varieties  you are considering.  \end{itemize}
       
   \vv If you organize your parameters appropriately you can ``geometrize'' your classification problem by recasting it as the problem of finding $K$-rational points on a specific  algebraic variety. 
   \vv

        In more technical vocabulary: you've framed a {\it representable moduli problem}---and the algebraic variety in question is called the {\it moduli space representing that   moduli problem}. \vv
  
        {\bf  Some classical examples---modular curves:} Fixing $N$ a positive integer and sticking to elliptic curves, the moduli spaces for rational torsion points or cyclic isogenies are  smooth curves defined over $\Q$:
 \[\begin{tikzcd}
	{\textrm{torsion points of order}\; N:} & {Y_1(N)} & {X_1(N)} \\
	{\textrm{cyclic isogenies of degree}\; N:} & {Y_0(N)} & {X_0(N)}
	\arrow[hook, from=1-2, to=1-3]
	\arrow[from=1-2, to=2-2]
	\arrow[from=1-3, to=2-3]
	\arrow[hook, from=2-2, to=2-3]
\end{tikzcd}\]        
   \vv

         The elliptic curves defined over $K$ {\it  possessing  a $K$-rational point of order $N$} are classified by the $K$-rational points of the affine curve $Y_1(N)$---and $X_1(N)$ is the  smooth projective completion of $Y_1(N)$  obtained by the adjunction of a finite set of cusps.  
         \vv And similarly, the classification of elliptic curves defined over $K$ {\it possessing  a $K$-rational cyclic isogeny of degree $N$} is related to the $K$-rational points of the affine curve $Y_0(N)$, which is a coarse moduli space.  The curve $X_0(N)$ is the smooth projective completion of the curve $Y_0(N)$.
         
         \vv
        
  The geometric formulation comes with a number of side-benefits.  Here are two:\vv
       \begin{enumerate}\item   If, say, the curve $X_0(N)$ is of genus $0$---noting that one of the cusps (${\bf \infty}$)  is defined over ${\Q}$, it follows that there is  a rational parametrization of that curve over ${\Q}$ which gives us a systematic account  (and parametrization);  that is, a $K$-rational parametrization of cyclic $N$-isogenies of elliptic curves---for any $K$.\vv
        \item  If it is of genus greater than 0,  one has a ${\Q}$-rational embedding (sending the cusp ${\bf \infty}$ to the origin) 
       $$X_0(N) \longhookrightarrow\ J_0(N)$$
       \noindent of the curve in its Jacobian, which allows us to relate questions about $K$-rational cyclic $N$-isogenies to questions about the Mordell--Weil group (of $K$-rational points) of the abelian variety $J_0(N)$.\vv
       
      \end{enumerate}      Besides  being able to apply all these resources of Diophantine techniques, there are the simple constructions that are easy to take advantage of.
    \vv
   For example, if you have  a moduli space  ${\mathcal M}$ whose $K$-rational points for every number field $K$ provides a classification of your problem over $K$, then, say, for any prime $p$ the set of  $K$-rational points of the algebraic variety that is the $p$-th symmetric power of   ${\mathcal M}$---denoted $\Symm^p({\mathcal M})$---{\it essentially} classifies the same problem  ranging over {\it all} extensions of $K$ of degree  $p$.  Given a variety  $V$ over a field $K$, by  a {\bf degree $n$ point (of  $V$ over $K$)} we mean a rational point of $V$ over some field extension of $K$ of degree $n$.   The degree 2 points are known as the quadratic points.\vv
   
    As an illustration of this, consider cyclic isogenies  of degree $N$ and noting that the natural ${\Q}$-rational mapping   
   \begin{equation}\label{Sym}  \Symm^p(X_0(N))\ \ \ \ \ \longrightarrow \ \ \ \ \ J_0(N)\end{equation}
   
    \noindent defined  by
    $$(x_1,x_2,\dots, x_p)\ \ \ \ \  \mapsto\ \ \ \ \     {\rm Divisor \ class\ of\ }\ \ \  [\sum_ix_i \ -\ p\cdot{\bf \infty}]$$
    \noindent has linear spaces as fibers, we get that the  classification problem of all cyclic $N$-isogenies of elliptic curves over all number fields of degree $p$ is geometrically related, again,  to the Mordell--Weil group of $J_0(N)$ over ${\Q}$.   
    \vv
        
A particularly nice example of this strategy carried out in the case of the symmetric square $\Symm^2$ of \emph{Bring's curve} is  in the appendix by Netan Dogra. {Bring's curve} is the smooth projective genus 4 curve in $\mathbb{P}^4$ defined  as the locus of common zeros of the following system of equations:
\begin{equation}\label{brings}
\begin{cases}
x_1 + x_2 + x_3 + x_4 + x_5 &= 0 \\
x_1^2 + x_2^2+ x_3^2 +  x_4^2 + x_5^2 &= 0 \\
x_1^3 + x_2^3+ x_3^3 +  x_4^3 + x_5^3  &= 0.  \end{cases}
\end{equation} It has no real points and thus no rational points. However, there are a number of points defined over $\Q(i)$, such as $(1:  i: -1: - i: 0)$. A natural question is thus if one can find \emph{all quadratic points} on Bring's curve. Dogra proves that all quadratic points are defined over $\Q(i)$ and produces the complete list of $\Q(i)$-rational points.

    \vv
    
 Part \ref{T} below concentrates on Ogg's torsion conjectures and the results that have emerged that are relevant to them.  In Part \ref{B} we review the broad uniformity conjectures (and results)  that have evolved from that work. Part \ref{CCK}  is a  discussion of the more recent method of Chabauty--Coleman--Kim  designed to compute rational points on curves by $p$-adic considerations; we focus specifically on the results achieved by this method for computation  of rational points on specific families of modular curves. 

    \part{Torsion and Isogenies}\label{T}
        \section{Ogg's Torsion Conjectures  (1973)}\label{etc}

         Ogg's Torsion Conjectures taken in broad terms can be formulated in terms of ``the geometrizations," as just described---i.e., in terms of ${\Q}$-rational points of modular curves---and the Mordell--Weil groups of abelian  varieties (i.e., of their Jacobians): \vv
       \subsection{${\Q}$-Rational torsion}\label{qrt}\
        \vv  
         \begin{quote}  {\bf Conjecture 1  (Ogg):}  An isomorphism class $\{C\}$ of finite groups occurs  as the torsion subgroup of the Mordell--Weil group of some  elliptic curve (defined over ${\Q}$)  {\it if and only if} the modular curve that classifies this problem is of genus zero{\footnote{ A form of this conjecture was  made by Beppo Levi in his 1908 ICM address in Rome. See \cite{17} which gives a wonderful account of the  story of Levi's engagement with (and his important results about)  the arithmetic of elliptic curves---all this being even before Mordell proved that the group of rational points of an elliptic curve over ${\Q}$ is finitely generated. Levi considers the tactic of  producing multiples of a rational point on an elliptic curve $\{n\cdot P\}\ n=1,2,3,\dots$ a ``failure" if it loops finitely---i.e., if P is a torsion point;  his aim is to classify such  ``failures." \vv}} .\end{quote}
       \vv

         Put in another way: an isomorphism class occurs  if and only if  it is expected to occur; i.e.,  if it  necessarily occurs, as a consequence of the ambient  geometry---this view being a continuing  guiding inspiration for number theory.\vv
        By ``geometry''  one means   the (algebraic) geometry of the curve $X_0(N)$.  For example, Ogg's article  \cite{14}  discusses the curious case of $X_0(37)$ which has two  noncuspical ${\Q}$-rational points, these being the images of the hyperelliptic involution (a non-modular involution) applied to the two cusps, both cusps being ${\Q}$-rational{\footnote{ See Section \ref{Ex} below.}}.  Ogg comments:\vv
         \begin{quote}  As Mazur and I are inclining to the opinion that $Y_0(N)$ has no ${\Q}$-rational  points except for a finite number of values of $N$, we are certainly interested in knowing when this sort of thing is going on, and in putting a stop to it if at all possible.\end{quote}
        \vv
        \subsection{ $\Q$-rational cyclic isogenies}\label{qri}
        \vv
        
  There are two different proofs of Conjecture 1. A major step in one of these proofs of Conjecture 1 is the full classification of  $\Q$-rational cyclic isogenies of prime degree; this is proved in   \cite{M2}:
  
\begin{thm}\label{1.1} Let $N$ be a prime number such that some elliptic curve over ${\Q}$ admits a
${\Q}$-rational $N$-isogeny. Then 

$$N=2, 3, 5, 7, 13 \ ({\rm the\  genus \ zero \ cases})$$
\centerline{ or}
$$N=11, 17, 19, 37, 43, 67,\ {\rm or}\  163.$$\end{thm}

This result was followed by a  sequence of papers of M.A. Kenku (\cite{19, 20, 21, 22})  that extends the  classification  to cyclic isogenies of any degree: 
\begin{thm}\label{1.2}The ${\Q}$-rational cyclic isogenies of degree $N$ of elliptic curves defined over ${\Q}$ only occur---and do occur---if  $1 \le N \le 19$ or if $N=21, 25, 27, 37, 43, 67,$ or $ 163$.  \end{thm} \vv  Following in the spirit of Ogg's original view of torsion points, {\it all} of these $N$-isogenies have ``geometric reasons'' for existing; e.g., the $37$-isogenies come by applying the hyperelliptic involution (it is non-modular!) to the cusps of $X_0(37)$. \vskip20pt
   
   \subsection{Rational torsion points on the Jacobians of  modular curves}
    Noting that  the cusps of $X_0(N)$ map to torsion points  of  its Jacobian,  $J_0(N)$, denote by $C_0(N) \subset J_0(N)_{\rm tors}  \subset   J_0(N)$ the subgroup generated by those cusps.     
  \[\begin{tikzcd}
	{\textrm{Cusps in}\; X_0(N)} \\
	{C_0(N)} & {X_0(N)} \\
	{J_0(N)_{\textrm{tors}}} & {J_0(N)}
	\arrow[from=1-1, to=2-1]
	\arrow[from=2-1, to=3-1]
	\arrow[from=1-1, to=2-2]
	\arrow[from=2-2, to=3-2]
	\arrow[hook, from=3-1, to=3-2]
\end{tikzcd}\]

      We have another, seemingly quite different type of conjecture:
      \vv
          \begin{quote}  {\bf Conjecture 2:}\label{conj2} Let $N$ be a {\it prime number}. We have:
          $$C_0(N)=J_0(N)_{\rm tors}({\Q})  \subset   J_0(N)({\Q})$$
         
               \end{quote}  \vv

        Put in another way: there are no ``unexpected''  ${\Q}$-rational torsion points  in $J_0(N)$: they all come from cusps.
   \vv

 Conjectures 1 and 2  are known.
 For Conjecture 1, see    in \cite[Theorem  8]{M1} and also \cite[Theorem 2]{M2}.  For Conjecture 2, see  \cite[Theorem 1]{M2}.   (Also see the broad survey of rational torsion results in  Andrew Sutherland's  \cite{16}).   That these conjectures  are interlinked is a long story,  as we discuss in Section \ref{rel}.
     \vv
 
   \subsection{ Conjecture 1}\label{conj1} Letting $C_n$ denote the cyclic group of order $n$,  the complete list of possible (isomorphism classes of)  finite groups that occur as torsion subgroups  of the  Mordell--Weil group  of ${\Q}$-rational points of elliptic curves are\vv   \begin{itemize} \item $C_n$ with $1\le n \le10$, and also $C_{12}$, and
 \item  the direct sum of $C_2$ with $C_{2m}$,  for $1\le m \le 4$.
      \end{itemize} \vv All these torsion groups occur infinitely often over $\Q$, since the corresponding modular curves are all genus zero curves   possessing a rational point.\footnote{See \cite{Tr} where it is proved that each of these groups appears as a possible torsion group over any quadratic field.}
      
       \begin{rem}\label{rem} Thanks to the work of Lo{\"i}c  Merel  \cite{Mer}, Joseph Oesterl{\' e},  Pierre Parent \cite{P} and others, we  have neat explicit upper bounds for the order of torsion points on elliptic curves   over number fields  of degree $d$.  For surveys of this work,  see \cite{D} and \cite{16}.\end{rem}
       \vv
       
  \section{The connection between algebraic torsion on elliptic curves and  rational  torsion on abelian varieties related to elliptic curves}\label{rel}
  
  \vv 
    The easiest way to explain this is to follow the ideas of the proof of Conjecture 1 in \cite{M2}, rather than the ideas in the earlier and quite different proof  in \cite{M0}.

            \vv
  
To set things up, let $N$ be a prime number such that $X_0(N)$ is of genus greater than $0$ and  let $J_{/ {\mathbb Z}}$  be the N{\'e}ron model of the Jacobian of $X_0(N)$ over ${\Q}$.  Let $$X_0(N)^{smooth}_{/ {\mathbb Z}} \ \ \stackrel{\ \ \iota}{\longhookrightarrow}\ \  J_{/ {\mathbb Z}}$$ be the smooth locus of the Zariski closure of    $X_0(N)_{\Q}$  in  $J_{/ {\mathbb Z}}$,  the embedding   being defined by  sending  the cusp ${{\bf \infty}}$---viewed as    ${\mathbb Z}$-valued section ${\bf e} \in X_0(N)_{/ {\mathbb Z}} $---to the `origin section' of    $J_{/ {\mathbb Z}}$.\vv

 An elliptic curve  $E$ with a cyclic isogeny of degree $N$ over  ${\mathbb Q}$ is represented by  a {\it noncuspidal} Spec ${\mathbb Z}$-valued section,  ${\bf x}$, of  $X_0(N)^{smooth}_{/ {\mathbb Z}}$ and hence (via $\iota$) also of the N{\'e}ron model $J_{/ {\mathbb Z}}$ of the Jacobian of $X_0(N)$ over ${\Q}$. \vv Suppose that such a rational point  ${\bf x}$ exists, and  denote by ${\bar {\bf x}}$ its image under the Atkin--Lehner involution $w_n:  X_0(N)\to X_0(N)$, the involution that exchanges the two cusp sections ${\bf 0}$ and  ${{\bf \infty}}$ of   $X_0(N)$.  Neither    ${\bf x}$ nor  ${\bar {\bf x}}$ are cuspidal sections of   $X_0(N)^{smooth}_{/ {\mathbb Z}}$.
\vv
 The articles \cite{M1} and \cite{M2} construct and discuss a specific smooth group scheme $A_{/ {\mathbb Z}}$ that is an optimal quotient{\footnote{ The group scheme $A$ is the relevant {\it Eisenstein quotient}---cf. loc.cit.  The term ``optimal quotient'' means that the kernel of  $J_{/ {\mathbb Z}}\to  A_{/ {\mathbb Z}}$ is a connected group scheme.}} of  $J_{/ {\mathbb Z}}$ for which these two properties are proven:\vv
 \noindent {\bf(1)}\   The generic fiber $A_{/ {\mathbb Q}}$  of  $A_{/ {\mathbb Z}}$ is an  abelian variety with {\it finite} Mordell--Weil group---i.e., $A({\mathbb Q})$   consists of rational torsion--and hence the image under $f$  of any   $ {\mathbb Z}$-valued section of $ X_0(N)^{smooth}_{/ {\mathbb Z} }$ is either trivial, or else generates a cyclic finite flat subgroup of  $A_{/ {\mathbb Z}}$;\vv
 \centerline{and:}\vv
\noindent {\bf(2)}\ The following diagram 

\[\begin{tikzcd}
	{\textbf{x}} & {X_0(N)_{/\mathbb{Z}}^{smooth}} & {J_{/\mathbb{Z}}} \\
	{\textbf{e}} & {X_0(N)_{/\mathbb{Z}}^{smooth}} & {A_{/\mathbb{Z}}}
	\arrow[from=1-1, to=1-2]
	\arrow[hook, from=1-2, to=1-3]
	\arrow[from=1-1, to=2-2]
	\arrow[from=2-1, to=2-2]
	\arrow["f", from=1-2, to=2-3]
	\arrow[from=1-3, to=2-3]
	\arrow["f", from=2-2, to=2-3]
\end{tikzcd}\]
\vv
has the property that\vv

  \begin{enumerate} \item[(i)]  the mapping   $f:X_0(N)^{smooth}_{/ {\mathbb Z} }\longrightarrow A_{/ {\mathbb Z}}$
is a formal immersion along the  cuspidal section ${\bf e}$, and
\item[(ii)] the diagram\vv\hskip70pt  $\xymatrix{ X_0(N)\ar[r]^f \ar[d]^{w_N}  & A\ar[d]^{-1}\\
 X_0(N)\ar[r]^f   & A }$ \vv 
 \noindent commutes, where ``$-1$" denotes the involution $z \mapsto  z^{-1}$.\end{enumerate} \vv So, by (i),  $f$ is a formal immersion along both cusp sections.   It follows that the image $\alpha\coloneqq f({\bf x}) \in A({\mathbb Z})$  is either\vv
 \begin{itemize} \item  the section defining the origin in the    group scheme $A_{/ {\mathbb Z}}$ or else \vv\item  is the generating section of a (nontrivial) cyclic finite flat subgroup scheme over ${\mathbb Z}$, call it:  ${\mathcal G}_{/ {\mathbb Z}} \subset A_{/ {\mathbb Z}}.$ 
 \vv  
 \noindent  By the classification of such group schemes we have that either      ${\mathcal G}_{/ {\mathbb Z}}$ is a constant (nontrivial ) group scheme, or else ${\mathcal G}_{/ {\mathbb Z}} \simeq \mu_2$   ($\mu_2 \subset {\mathbb G}_m$ being the kernel of multiplication by $2$ in the multiplicative group scheme  ${\mathbb G}_m$). \vv
\end{itemize}   These possibilities also hold, of course, for the ``conjugate section''  ${\bar \alpha}\coloneqq f({\bar{\bf x}}) \in A({\mathbb Z})$: it is either the trivial section or it generates a finite  flat group scheme   ${\bar {\mathcal G}}_{/ {\mathbb Z}} \subset A_{/ {\mathbb Z}}$  that is either a constant group scheme or $\mu_2$.
 \vv
\begin{lem} Neither ${\alpha}$ nor ${\bar \alpha}$ are the trivial section of $A_{/ {\mathbb Z}}$.\end{lem}

\begin{proof} Since $f$ is a formal immersion along the cuspidal sections  if  ${\alpha}$ or ${\bar \alpha}$ were the trivial section  we would be led to a contradiction, as illustrated by the following diagram (taken from \cite[p. 145]{M2}):

\vv
     \begin{center}
 \begin{figure}[!h]
 \caption{Contradiction:  {\it  intersection with the trivial section}}
  \vskip10pt \includegraphics[width= 0.8\textwidth]{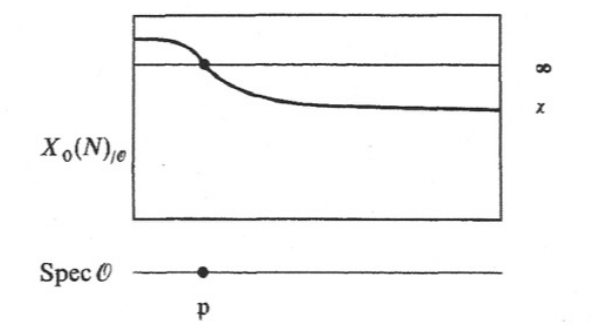}\end{figure}\end{center}
   
   \noindent in that the image of the two depicted sections would converge onto the origin section of $A$ contradicting formal smoothness along ${\bf e}$.\end{proof} 

 So  ${\alpha}$ and ${\bar \alpha}$ are generators of nontrivial  group schemes     $ {\mathcal G}$ and    ${\bar {\mathcal G}}$  respectively, these being either constant or $\mu_2$. 
 \begin{itemize} \item  If   $ {\mathcal G}$ and    ${\bar {\mathcal G}}$   are constant group schemes, then    ${\alpha}$ and ${\bar \alpha}$ are sections of $A$  over $\Z$   disjoint (as schemes) from the trivial section of $A$ and therefore  ${\bf x}$ and ${\bar {\bf x}}$ are disjoint (as schemes)  from the cuspidal sections of $X_0(N)_{/ {\mathbb Z}}$.    It follows that the elliptic curves $E$  and ${\bar E}$ that are classified by  ${\bf x}$ and ${\bar {\bf x}}$  have {\it potentially good reduction everywhere.} \vv
 \item  And if   ${\alpha}$ or ${\bar \alpha}$ generates  a subgroup isomorphic to $\mu_2$,  since $\mu_2$ is {\'e}tale outside the prime $2$ it follows that   $E$  or  ${\bar E}$  would  have {\it potentially good reduction except for the prime $2$.}

\end{itemize}   
    
          \vv  Even though this is just the start of a sketch of the proof of Conjecture 1 in \cite{M2},  from what we've just  described, we can prove
          
          \vv
          \begin{thm}   The only prime numbers $N$ for which there exist elliptic curves over ${\mathbb Q}$ with rational torsion  points of order $N$ are: $$N=2,3,5,7.$$\end{thm}\vv
          \begin{proof}  
          First note that   $X_1(N)$ is of genus $0$ for  $N=2,3,5,7$ so there are infinitely many elliptic curves with rational torsion points of order $N$ for  these primes. That list of primes  {\it and}\ $13$ are precisely  the primes for which $X_0(N)$  is of genus $0$. Since  the curious prime $N=13$ is taken care of by \cite{MT} where it is proven that there are no rational points of order $13$ on elliptic curves over ${\Q}$, to prove the theorem we may suppose  $N$ to be different from   $N=2,3,5,7, 13$; equivalently,  that $X_0(N)$ is of genus greater than $0$ so the discussion above applies. \vv
          
          In particular,  we assume that $E$ is an elliptic curve over ${\Q}$, of potential good reduction away from $p=2$, and possessing a rational point of order  $N=11$ or $N \ge17$, where $N$ is a prime. Since it has  such a rational point, 
 the N{\'e}ron model of $E$ over ${\mathbb Z}$   contains a constant subgroup scheme ${\mathcal Z}$ isomorphic to $\Z/(N\cdot \Z)$.
          
          \vv  For $p$ a prime, let $E_p$ denote the  fiber at the prime $p$ of the  N{\'e}ron  model of $E$, so $E_p$ is a group (scheme) of finite order over the prime field ${\mathbb F}_p$.   Since the specialization of  ${\mathcal Z}$  to $E_p$ defines $N$ distinct   ${\mathbb F}_p$-rational points of $E_p$   it follows that \begin{equation}\label{lb}  N\ {\rm divides}\  |E_p({\mathbb F}_p)|.\end{equation}
 \vv         
          If $p>2$, since $E$ is of potentially good reduction, in the terminology of the  theorem of Kodaira and N{\'e}ron (cf. Theorem 15.2  and Table 15.1 in  Section 15 in Appendix C of \cite{Si}) we have that  $E_p$ is not of  multiplicative type---i.e., of  type I$_\nu$  or  I$_\nu^*$  for any $\nu > 0$.   So either:\begin{itemize}\item $p$ is a prime  of good reduction for $E$, or \item it is of additive reduction.\end{itemize} 
          \vv   If $E$ has additive  reduction at $p$   (i.e., the N{\'e}ron fiber at $p$ is of  one of the types II, III, IV or I$_0^*$, II$^*$, III$^*$, or IV$^*$; see Table 15.1 in loc.cit.)    then $E_p$  is an extension of the additive group ${\mathbb G}_a$  over ${\mathbb F}_p$ by a finite group of order $\le 4$. In particular $|E_p({\mathbb F}_p)|$ is divisible by $p$ and is $\le 4p$. It follows that Equation \ref{lb}, applied to  the prime $p=3$ already shows that $E$ cannot have additive reduction  at $p=3$ for the primes $N$ we are considering, so it must have good reduction---i.e., be an elliptic curve---at $p=3$.  But since any elliptic curve over  ${\mathbb F}_3$ has at most seven  ${\mathbb F}_3$-rational points, we see, by \eqref{lb} that $N$ is either $2$, $3$, $5$, or $7$.     \end{proof}
          
 \vv A significantly more detailed outline of the proof  of Conjecture 1 is given as Steps 1--4  on pages 132, 133 of  \cite{M2}---the {\it full proof} itself being in the body of that paper.
 \part{Boundedness of torsion and isogenies  in more general contexts}\label{B} \vv 
 
              Conjecture 1 having been completely resolved in the case of elliptic curves  has  inspired  more general uniform boundedness expectations for rational points; e.g.,  for abelian varieties $A$ over number fields $K$:  conjectures that the order of the torsion group of an abelian variety over a number field can be bounded in terms of the dimension of the variety and the number field; and still stronger versions:  that the torsion is bounded in terms of the dimension of the variety and the degree of the number field.
              \vv
             
            Moreover,  it is striking how few additional  isomorphism classes of $K$-rational torsion subgroups of elliptic curves can occur in elliptic curves over quadratic  and cubic number fields $K$:
       
     \section{Torsion on elliptic curves over quadratic and  cubic number fields}  
      \begin{thm}[Momose, Kenku, Kamienny, {\cite{K1, K2, 9, 10, K3, K4, K5, K6, 22, K7, 12.5, Mo1, Mo2}}] \label{MKK}  Let $K$ range through all quadratic  number fields,  and $E$  all elliptic curves over  these fields.
Then the torsion subgroup $E(K)_{\rm tors}$  of $E(K)$ is isomorphic to one of the
following 26 groups:
\begin {itemize} \item $C_n$ for $1\le n \le 18, \ n \ne 17,$
 \item  the direct sum of $C_2$ with $C_{2m}$  for $1\le m\le 6$,
\item  the direct sum of $C_3$ with $C_{3m}$ for $m=1,2$,
\item $C_4\oplus C_4$.
\end{itemize} \end{thm}  \vv

     \begin{thm}[Derickx, Etropolski, Van Hoeij, Morrow, Zureick-Brown {\cite{DEHMZ}}]\label{DEHMZ-B} Let $K$ range through all cubic  number fields,  and $E$  all elliptic curves over  these fields. Then the torsion subgroup $E(K)_{\rm tors}$  of $E(K)$ is isomorphic to one of the
following 26 groups:
\begin{itemize} \item   $C_n$ for $1\le n \le 18, \ n \ne 17,$
 \item  the direct sum of $C_2$ with $C_{2m}$  for $1\le m\le 7$,
 \item $C_{20}, C_{21}.$
\end{itemize}
There exist infinitely many ${\Q}$-isomorphism classes for each such torsion subgroup except for $C_{21}$.
In this case, the base change of the elliptic curve with LMFDB label \href{https://www.lmfdb.org/EllipticCurve/Q/162/c/3/}{162.c3} to ${\Q}(\zeta_9)^+$
 is the unique elliptic curve over
a cubic field  $K$ with  $K$-rational torsion group isomorphic to  $C_{21}$. \end{thm}\vv

         \subsection{Conjecture 2 expanded and related to cuspidal subgroups}\
         
       \vv
       Recall:
       \begin{quote}  {\bf Conjecture 2:}\label{conj2.} Let $N$ be a {\it prime number}. We have:
          $$C_0(N)=J_0(N)_{\rm tors}({\Q}) \ \ \subset \ \   J_0(N)({\Q}).$$
         
               \end{quote}  \vv
        \begin{itemize}\item The order of the $C_0(N)$  had been computed for square-free $N$ thanks to  Kubert and Lang \cite{Kubert-Lang}, and Takagi \cite{takagi}.   In this case  (i.e., $N$ square-free) the set of cusps are ${\Q}$-rational. \vv

 \item Ohta \cite{O1, O2} has proved a generalization of this conjecture in the context of square-free $N$.
That is, he proved that the $p$-primary parts of $J_0(N)_{\rm tors}({\Q}) $ and of $C_0(N)$
are equal for $p \ge 5$   and $p=3$ if $3$ doesn't divide $N$.\vv    

Related to this, see \cite{L}, \cite{CES}, \cite{Y}, and \cite{R1},  \cite{R2}.   And very recently the PNAS article \cite{RW}  ({\it Another look at rational torsion of modular Jacobians}) by Ken Ribet and Preston Wake appeared, giving another approach to this issue.\vv
      \item In the more general context of $N$ not square-free, the cuspidal subgroup
of $J_0(N)$ may not consist entirely of rational points; nevertheless:

\vv 
{\bf Conjecture $2^*$:}  $$J_0(N)_{\rm tors}({\Q})\  =C_0(N)({\Q})\subset C_0(N).$$\end{itemize}

\vv
          \subsection{Conjecture 2 further expanded}\
          \vv
          
        Let   $X$  be a modular curve  (over $\mathbb{Q}$)  and    ${\mathcal J}$    the Jacobian of  $X$. Let  $${\mathcal C} \subset {\mathcal J}$$
be the finite {\'e}tale subgroup scheme of  ${\mathcal J}$  generated by the cusps. Let $K/{\mathbb{Q}}$ be the field cut out by the action of Galois on  ${\mathcal C}$.  Thus there's an exact sequence
  $$0\to  \Gal({\overline{\mathbb{Q}}}/{K}) \to \Gal({\overline{\mathbb{Q}}}/{\mathbb{Q}})\to \Aut({\mathcal C}({\overline{\mathbb{Q}}})).$$
\vv
Define the {\bf cuspidal defect} of  $X$ to be the cokernel of
\begin{equation}\label{cok}{\mathcal C}({\overline{\mathbb{Q}}})= {\mathcal C}(K)  \longhookrightarrow   {\mathcal J}(K)_{\rm tors}.\end{equation}
\vv
{\bf Conjecture $2^{**}$:}  Let $X$ be either  $X_0(N)$ or $X_1(N)$ for some $N \ge 1$.   The ``cuspidal defect'' of $X$ is trivial.
\vv  

\vv

          \section{Remarkable ``Diophantine Stability''} 
\begin{defn}
\label{def0}
Let $L/K$ be an extension of (number) fields, and $V$ an algebraic variety 
defined over $K$. Denote by $V(K)$ the set of $K$-rational points of $V$.
Say that $V$ is {\bf Diophantine Stable} for $L/K$, or 
$L/K$ is {\bf Diophantine Stable} for $V$,
if the inclusion $V(K) \hookrightarrow V(L)$ is an isomorphism, 
i.e.:
{\it if $V$ acquires no new rational points after passing from $K$ to $L$.}
\end{defn}
     \vv
     Note that Theorem \ref{MKK} tells us that:
         \vv
         \begin{cor}\label{ds1} For all but finitely many positive numbers $N$, the curve $Y_1(N)$   (over $\Q$) is Diophantine Stable for {\bf all} quadratic extensions $L/{\Q}$.\end{cor}
     
\vv
       This is striking and suggests that Diophantine Stability is a common feature.{\footnote{ Filip Najman suggested that one might add a comment that the Diophantine Stability phenomenon of Corollary 
       \ref{ds1} holds more generally over number fields of any degree, given the results referred to in Remark \ref{rem} in Subsection \ref{conj1} above.\vskip10pt}}
       
       \vv Consider:
       
\begin{thm}[Theorem 1.2 of \cite{MR1}]
\label{intro:1}
Suppose $A$ is a simple abelian variety over $K$ and all $\overline{K}$-endomorphisms of 
$A$ are defined over $K$. Then there is a set ${\mathcal S}$ of rational primes with
positive density such that for every $\ell \in {\mathcal S}$ and every $n \ge 1$, there 
are infinitely many cyclic extensions $L/K$ of degree $\ell^n$ such that 
$A(L) = A(K)$.

\vv

If $A$ is an elliptic curve without complex multiplication, then ${\mathcal S}$ can be 
taken to contain all but finitely many rational primes.
\end{thm}

\vv
\noindent  and this  is surely not the last word regarding the extent of Diophantine Stability, specifically if the base field $K$ is ${\Q}$ and  if $A=E$, an elliptic curve over ${\Q}$.   We conjecture that any such $E$ is Diophantine stable for all but finitely many Galois extensions of prime degree greater than $5$.

    \vv   
          
\section{{\it ``Expected'' and ``Unexpected''} $L$-rational cyclic isogenies for $L$ ranging through quadratic fields}  
    \vv

      What about uniformity results regarding cyclic $N$-isogenies of elliptic curves ranging over {\it all} quadratic fields?  This  question has been addressed in \cite{6} and generalized to arbitrary number fields in \cite{BM}. \vv
      
        One source of uniformity theorems consists of consequences of a general theorem of Faltings: 
         \begin{thm}\label{Falt} Let $A$ be an abelian variety defined over a number field $L$. Any
closed subvariety of $A$ defined over $L$ which is the Zariski-closure of its set of
$L$-rational points is a finite union of translates of
abelian subvarieties of A.
\end{thm}
\vv
This was first proved by Faltings \cite{F94}  following \cite{F91}; also  see  \cite{Mc}.
Techniques of \cite{F91}\ have, as their starting point, Vojta's proof \cite{V91}
of the classical Mordell Conjecture. See also \cite{V93}.
For further discussion of this in the context of generalization(s) of the classical Mordell Conjecture---with references listing the people who also worked on this, see  \cite{M3}. \vv
  A corollary of a theorem of Faltings  combined with  \cite{HS}  is:
                \begin{cor}[Faltings, Harris--Silverman\footnote{See also the discussion and Theorem 1.3 of \cite{KV}, which handles the case when $X$ is a double
cover of an elliptic curve of positive rank.}] Let $K$  be a number field and $X$ a curve defined over  $K$.
        Then $X$ is Diophantine Stable for all but finitely many quadratic extensions $L/K$ unless $X$ is
        \begin{itemize}\item of genus $0$ or $1$, or
        \item hyperelliptic or bielliptic  (over $K$).\end{itemize}\end{cor} 
        \vv
        
        \begin{proof}  Assume that $X$ is not of genus $0$ or $1$, or
      hyperelliptic or bielliptic.  So the genus of $X$ is $\ge 3$. Let ${\rm Div}^0(X)$ denote the group of divisors of degree zero in $X$; let $A:={\rm Pic}^0(X)$  and consider the natural diagram\vv
      
     \hskip90pt  $\xymatrix{ \Symm^2(X)\ar[r]^\iota\ar[rd]^j &{\Div}^0(X)\ar[d]^\pi\\
      \ &  A}$
      
      \vv where \vv\begin{itemize} \item  for $\{x_1,x_2\} \in \Symm^2(X)$ define  $$\iota: \{x_1,x_2\}\ \  \mapsto\ \  [x_1]+[x_2] - 2\cdot[0]\ \  \in \ \ \Div^0(X),$$
      \and \vv \item for $D\in {\rm Div}^0(X)$ define $$\pi: D \mapsto \ {\rm its\ linear\ equivalence\ class\ in \ }{\rm Pic}^0(X)=A.$$\end{itemize}\vv
       Noting that the connected components of the fibers of $\pi$ and of $j$  are (essentially by definition) rational varieties, and that since $X$ is assumed to be non-hyperelliptic $\Symm^2(X)$ contains no rational curve, it follows that we have an injection $j:\Symm^2(X) \hookrightarrow A$ .  Let $W\subset   A$ denote the image of  $\Symm^2(X)$ in $A$ (under the injection  $j$).  Since $X$ is of genus $\ge 3$,  $W$ is a proper subsurface in $A$. If it contained an elliptic curve, then $X$ would be bielliptic \cite{HS}, but since that's not the case, $W$ contains no (positive-dimensional!) abelian subvarieties---so by Faltings theorem has only finitely many $K$-rational points, establishing the corollary. 
(Compare Theorem 1.3 of \cite{KV}.)\end{proof}

\begin{figure}[!h]
\caption{Types of rational points on $\Symm^2(X_0(N)$---ignoring the possible existence of abelian surfaces}\label{fig1}
\centering \includegraphics[width=0.75\textwidth]{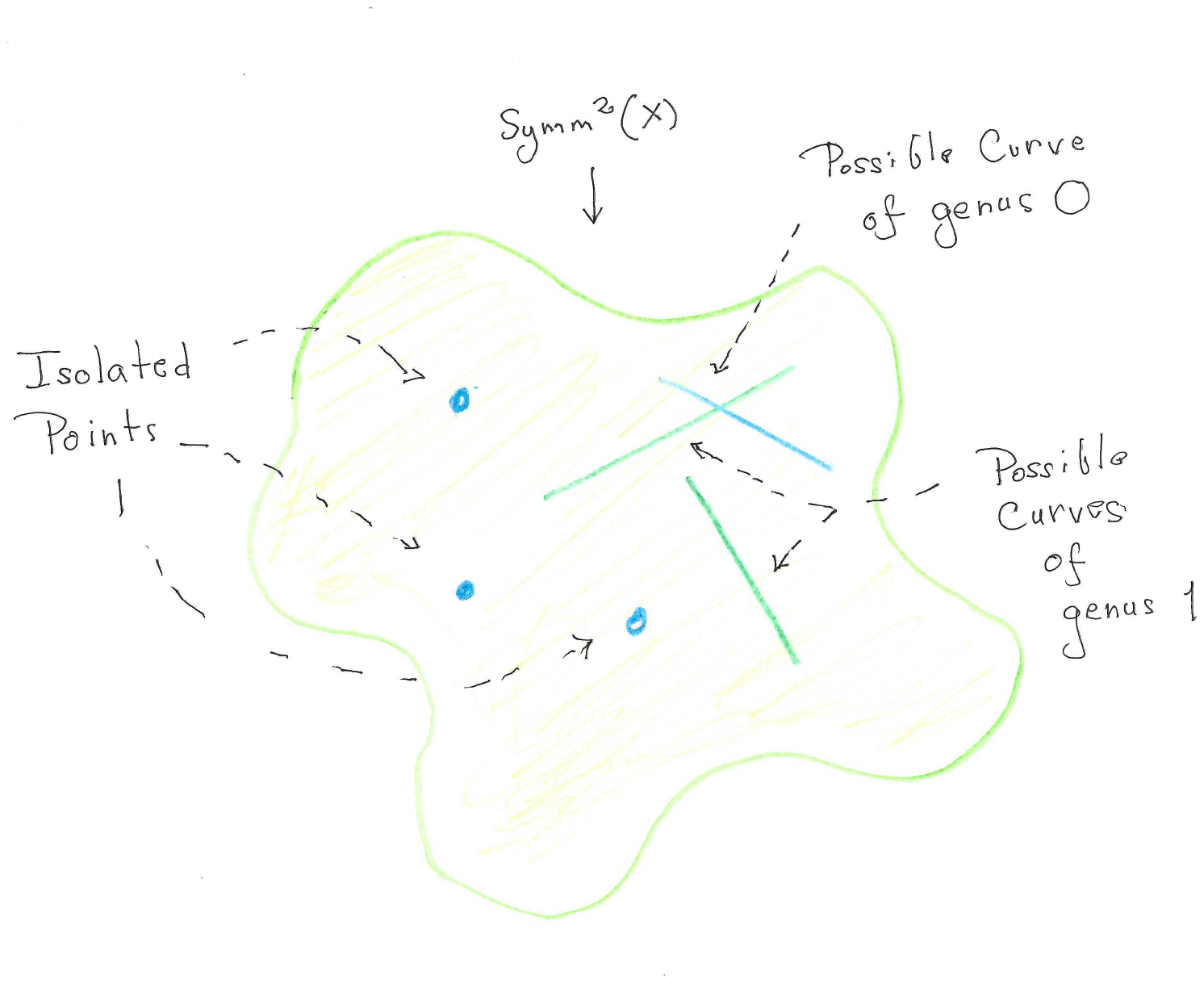}
\end{figure}
  
    \subsection{Isolated quadratic points}
Call the set of quadratic points of $X$ that are not among (infinite) systems of parametrized quadratic points {\bf isolated points}. (See \cite[Definition 4.1] {BELOV} and \cite[Definition 3.1]{VV} for the precise definitions and for important relevant discussion. See also \cite{HS, KV}.)   

Infinite systems of parametrized quadratic points deserve to be called ``expected  quadratic points (over $K$)  in $X$'' given the geometry of the situation.

\vv

But when $X=X_0(N)$ for some $N$ and $K ={\Q}$ there  may also be a few other points  of $X_0(N)$ over quadratic imaginary fields $ {\Q}({\sqrt{d}})$ of class number $1$; i.e.,  
$$d \ =\ -1, -2, -3, -7, -11, -19, -43, -67, -163$$

\noindent that   deserve the title ``expected.''   Namely,  if $E$ is an  elliptic curve over ${\Q}$ that is CM  with CM field   $K \coloneqq {\Q}({\sqrt{d}})$   (with $d$ in the above list) then for any positive integer $N$ with the property that all of its prime divisors are  (unramified and) split in $K$, $E$ has a $K$-rational cyclic isogeny of degree $N$; hence is classified by a $K$-rational point of $X_0(N)$.  Such a point is therefore also ``expected.''    So: \vv 
       \subsection{Sporadic quadratic points}
\begin{defn}  Call a quadratic point of $X_0(N)$  {\bf sporadic (quadratic)}{\footnote{Note that this differs from the concept  {\it sporadic} given in Definition 4.1 (4) of    \cite{BELOV}.}} if:

\begin{itemize} \item it is not a cusp, \newline \centerline{and}\vv
\item{\it is}  isolated; in particular,\vv\begin{itemize} \item  is {\it not}  the inverse image of  a  ${\Q}$-rational  point  in ${\mathbb P}^1$  via a hyperelliptic covering  (i.e., a degree $2$ mapping $X_0(N)\to {\mathbb P}^1$),  in the case where  $X_0(N)$ is hyperelliptic,  \vv  \centerline{and}  \vv \end{itemize}  \item is not  a point of $X_0(N)$ classifying a CM elliptic curve and cyclic isogeny of degree $N$ as described above.\end{itemize}  \end{defn}

See Figure \ref{fig1} for an illustration of these points.
   
   \begin{conj} Ranging over all curves $X_0(N)$ for $N \in {\mathbb Z}_{\ge 1}$  there are only finitely many sporadic quadratic points.\end{conj} \vv

    Surely  all of us agree with the spirit of  the quotation of Ogg's view regarding rational torsion  in Section \ref{etc}.  That is, we're interested ``in knowing when this  [sporadic quadratic points] sort of thing is going on, and in putting a stop to it if at all possible."
   \vv
   Thanks to the recent work of a number of people, {the sporadic quadratic points of all of the curves  $X_0(N)$ that are hyperelliptic or bielliptic have been computed}, {as we will discuss in the next section.

   \vv Sheldon Kamienny made the following comment:\vv
    \begin{quote} {\it The existence of sporadic points always left me scratching my head.  Do they fit into a framework, or is it just nature being unkind?} \end{quote} 
 \vv

   \section{Sporadic quadratic points on  hyperelliptic $X_0(N)$}
A classical theorem of Ogg \cite{14} gives the nineteen values of $N$ for which $X_0(N)$ is  hyperelliptic  (we take hyperelliptic to require that the genus is $>1$):\vskip10pt 
$ \begin{array}{ccccccccccc}
     N: & {\bf 22} & 23 & {\bf 26} &  {\bf 28} & 29 & {\bf 30} & 31& { \bf 33} & {\bf 35} & {\bf 37}\\ 
    {\rm  genus:} & 2 & 2 & 2& 2& 2 & 3 & 2& 3& 3& 2\end{array}$
    \vskip10pt 
$ \begin{array}{cccccccccc} 
N: &  {\bf 39} & {\bf 40} & 41 & 46 & 47& {\bf 48} & {\bf 50} & 59 & 71\\
   {\rm  genus:}    & 3 & 3 & 3 & 5 & 4 &3 & 2   & 5 & 6\end{array}$
\vskip10pt 

   The levels $N$ that appear in boldface above are those values of $N$ such that $X_0(N)$ is bielliptic as well as hyperelliptic. All sporadic  quadratic points for any of those modular curves $X_0(N)$  (except for $X_0(37)$)    have been computed by Peter Bruin and  Filip Najman  in their article \cite{8} (which has other interesting results as well).  The case of  $X_0(37)$ is taken care of in  Section 4 of Josha Box's paper  \cite{7}, in which all sporadic quadratic points have also been computed for the curves   $X_0(N)$ with  $N=43,  53, 61, 65$, these being bielliptic curves covering elliptic curves of positive Mordell--Weil rank.
   
   \begin{prop}[Bars \cite{6}]\label{list}  These are the values of $N$ for which  $X_0(N)$ is  of genus $>1$ and bielliptic (over ${\Q}$):
   \vskip10pt 
 $ \begin{array}{ccccccccccccccccc}  
  22 &26 & 28 & 30 & 33  & 34 & 35 & 37 & 38  & 39 & 40 & 42 &43 & 44 & 45 & 48 & 50\\
  51 & 53 & 54 & 55 & 56 & 60 & 61 & 62 & 63 & 64 &  65  & {69} & 72 & 75 &  {79} &  81 & 83 \\
  {89}& {92} & {94} & 95 &{101} & 119 & {131}    &\   &\ & \ &\   \end{array}$\end{prop}
    \vskip10pt 

    Until {\it very recently} there remained a dozen entries in the above table for which we did not  know the set of their isolated  quadratic points. Thanks to Filip Najman and Borna   Vukorepa \cite{NV}, we now have  computation of the isolated quadratic points  for all bielliptic curves $X_0(N)$  (as we also do for all hyperelliptic $X_0(N)$).

    \section{Exceptional quadratic points}
  \vv  
      Let $N$ be prime, and $w_N\colon X_0(N)\to X_0(N)$ the Atkin--Lehner involution.  This involution is given by sending  a pair  (representing a point in $X_0(N)$) $$(E, C_N\stackrel{\alpha}{\hookrightarrow} E)$$
      \vv\noindent---consisting of an elliptic curve $E$ and $C_N$  a cyclic subgroup of order  $N$---
      \vv\noindent to   the pair
      
      $$(E', C'_N\stackrel{\alpha'}{\hookrightarrow} E').$$\vv
      Here       
      $E' \coloneqq E/C_N$   and $C'_N \coloneqq E[N]/C_N$   (where $E[N]$ is the kernel of multiplication by $N$ in $E$).
\vv
Forming the quotient
      $$X_0(N)^+ \coloneqq X_0(N)/{\rm action\ of\ } w_N,$$
\noindent we get the double cover 
       $$X_0(N)  \stackrel{\pi}{\longrightarrow}  X_0(N)^+.$$
       
   \begin{defn}   For $N$ an integer where   $ X_0(N)^+$  of genus $>1$,    \begin{enumerate} \item call a ${\Q}$-rational  point  of $X_0(N)^+$  {\bf exceptional} if  it is neither a cusp nor a point classifying  a CM elliptic curve;  \item  call a quadratic  point $P$  of $X_0(N)$  {\bf exceptional} if it is not defined over ${\Q}$ (i.e., it is an honest quadratic point)  and the image of $P$ in   $X_0(N)^+$ is an exceptional ${\Q}$-rational  point.\end{enumerate}\vv \end{defn}
       \vv
  {\it Exceptional points }  deserve the adjective,    since they have the intriguing structure  of a duo of cyclic $N$-isogenies:

\[\begin{tikzcd}
E & {E'} & {\textrm{and}} & {E'} & E
	\arrow["N", leftrightarrow, from=1-1, to=1-2]
	\arrow["N", leftrightarrow, from=1-4, to=1-5]
\end{tikzcd}\]

    This structure can also be combined into a single abelian surface defined over ${\Q}$:
    $$A \coloneqq E\times E'$$ 
    \noindent  endowed with  an endomorphism:  $$``{\sqrt{N}}": (x,y)\mapsto (\alpha'(y), \alpha(x)).$$\vv
 
 What tools do we have to compute     the exceptional ${\Q}$-rational points on $X_0^+(N)$?
\vv
\part{The method of  Chabauty, Coleman and Kim}\label{CCK}\vv

The classical method of Chabauty--Coleman {(see Section \ref{chabauty-coleman} below)} computes  a usable bound for the number of rational points on a curve $X$  (of genus $>1$) provided  that the rank $r$ of the Mordell--Weil group of the Jacobian of $X$ is strictly less than its genus $g$.
\vv    But the Birch and Swinnerton-Dyer conjecture predicts that   (for $N$ prime) the rank $r_0(N)^+$ of {$J_0(N)^+(\Q)$}, the Mordell--Weil group of the Jacobian of $X_0(N)^+$, is greater than or equal to $g_0(N)^+$, the genus of {$X_0(N)^+$}.   So this classical method can't be brought to bear here. \vv

 Computationally, we have many examples where there's actual equality: $$r_0(N)^+=g_0(N)^+.$$   (Indeed, this is true for all $N < 5077$ for which $g_0(N)^+ > 1$.)  Happily, for exactly  such cases---i.e., for curves $X$ of genus $>1$ with $r=g$---we have the more recent   quadratic Chabauty--Coleman--Kim method that  offers 
a new approach to compute the set of all ${\mathbb Q}$-rational points{\footnote { We think it is reasonable to conjecture that the average value of the ratios  $$\frac{r_0(N)^+}{g_0(N)^+}$$ is $1$; e.g., as $N$ ranges through prime values;  are these ratios bounded?}}.  For example see \cite{4.2, 2, 3, 4} {(and Section \ref{qc})}. Indeed, there are also two new viewpoints on quadratic Chabauty, the geometric perspective of Edixhoven--Lido \cite{EL} and the ($p$-adic) Arakelov theoretic one of Besser--M\"uller--Srinivasan \cite{6.1}. We will say more about these in the following section.
   \vv
   
   The list of curves $X_0(N)^+$  of genus 2 or 3   with $N$ prime is a result of Ogg. 
 We have the following:
       \vv
    {\bf Theorem} (Ogg)   For  $N$\ prime,  $X_0(N)^+$ is of genus 2 if and only if $$N \in  \{67, 73, 103, 107, 167, 191\}$$ and it has genus 
3 if and only if
$$N \in \{97, 109, 113, 127, 139, 149, 151, 179, 239\}.$$ \vv

Elkies and Galbraith \cite{G1}  found {exceptional rational points} on  $X_0^+(N)$   for $N = 73, 91, 103, 191$ and $N = 137, 311$ {(which are of genus 4)}.  In  \cite{3}, it was shown that the only prime values of $N$ with  $X_0^+(N)$ of  genus 2 or 3 that have an {exceptional} rational point are $N = 73, 103, 191$ (all genus 2). In particular,  for prime $N$, if $X_0^+(N)$ is of genus 3,
it has no {exceptional} rational points.   Ad\v{z}aga, Arul, Beneish, Chen, Chidambaram, Keller, and Wen \cite{1.1} showed that the only prime values of $N$ with $X_0^+(N)$ of genus 4, 5, or 6 that have an {exceptional} rational point are $N = 137$ and 311.  {Thus for all of the above values of $N$, we have a complete understanding of the exceptional quadratic points on $X_0(N)$.} \vv

We briefly discuss the work of \cite{1.1} on the genus 4 curve $X_0(311)^+$.   Using the canonical embedding, a model for $X_0(311)^+$ is given by the following equations in $\mathbb{P}^3$:

\tiny
\begin{align*} &X^2 + W Y - 2 X Y + 2 Y^2 + 7 X Z - 8 Y Z + 13 Z^2 = 0, \\
&W X^2 - 2 W X Y + X^2 Y - W Y^2 - X Y^2 - 2 Y^3 + W^2 Z + 6 W X Z  - X^2 Z - W Y Z + 5 X Y Z + 4 Y^2 Z + 7 W Z^2 - 4 X Z^2 - 2 Z^3 = 0.\end{align*}
\normalsize
Using quadratic Chabauty (see Section \ref{qc}) at $p=5$ on a plane model, they show that there are precisely five rational points on the curve:

\begin{center}
\begin{tabular}{c|c}
\hline
rational point on  $X_0(311)^+$  & type of point \\
\hline
$(1 \colon 0 \colon 0 \colon 0)$&  cusp \\
$(1 \colon -1 \colon -1 \colon 0)$ & CM, $D = -11$ \\
$(1 \colon 2 \colon -1 \colon -1)$ & CM, $D=-19$ \\
$(2 \colon 0 \colon -1 \colon 0)$ & CM, $D=-43$ \\
$(6 \colon 8 \colon -1 \colon -2)$ & exceptional \\
\hline
\end{tabular}
\end{center}

\bigskip

Galbraith \cite{G1} had earlier computed that the $j$-invariant of the $\mathbb{Q}$-curve corresponding to the exceptional point is
\begin{align*}
    j = &\,31244183594433270730990985793058589729152601677824000000 \\
    &\quad \pm 1565810538998051715397339689492195035077551267840000\sqrt{39816211853}.
\end{align*}

\vv
See also the survey article  \cite{GLQ}  (and \cite{G1, G2}) in which exceptional points found by Elkies and Galbraith are defined and studied  in the context of  ${\mathbb Q}$-curves;  and for the list of the seven known exceptional $N$-isogenies, these being  rational over  a quadratic field  of discriminant $\Delta$:\vv

\begin{center}
  $ \begin{array}{c|c|c}
  \hline
  N   & g &    \Delta\\
  \hline
  73  &    2 &   -127    \\
  91 & 2 & -3 \cdot 29 \\
  103 &    2 &  5\cdot 557   \\
  125 & 2 & 509 \\
  137  &4    &  -31159     \\
  191  &    2 &  61\cdot 229\cdot 145757   \\
   311  &   4 &  11\cdot 17\cdot 9011\cdot 23629 \\
   
    \hline  \end{array}$ \end{center}

\vv
By the work of \cite{1.1, 1.2, 1.3, 2, 3} this gives a complete list of exceptional isogenies arising from rational points on the curves $X_0(N)^+$ of level $N$ and genus at most 6.  Are these the only exceptional  isogenies?  There's lots to be done.
\vv
\section{The method of Chabauty}
   The aim of this classical method is to prove finiteness of the set of ${\Q}$-rational points of  a curve $X$ of genus $g>1$ under the assumption that the rank  $r$ {of the Mordell--Weil group of the Jacobian $J$ of $X$} is small; specifically, if it is strictly less than $g$. \vv  One can assume that $X$ has at least one ${\Q}$-rational point, for otherwise the job is done. Choosing  a rational point $b \in X(\Q)$, form the Abel--Jacobi embedding 
\begin{align*}\ \ \ \ i_b: X &\rightarrow J\\ 
 P &\mapsto [(P) - (b)].\end{align*} 

For any prime $p$  viewing   $J({\Q}_p)$ as a $p$-adic analytic group (of dimension $g$) containing the  Mordell--Weil group $J({\Q})$ as subgroup, denote by $\Gamma_p \subset  J({\Q}_p)$ the topological closure of $J({\Q})$ in $J({\Q}_p)$  noting that, given our hypothesis,  its dimension is at most $r$, which is less than $g$. We have

\[\begin{tikzcd}
	{X(\mathbb{Q})} & {\Gamma_p} \\
	{X(\mathbb{Q}_p)} & {J(\mathbb{Q}_p)}
	\arrow["{i_b}", from=2-1, to=2-2]
	\arrow[hook, from=1-1, to=2-1]
	\arrow["{i_b}", from=1-1, to=1-2]
	\arrow[hook, from=1-2, to=2-2]
\end{tikzcd}\]

      $X({\Q}_p)$ is a (proper) $p$-adic analytic subvariety of   $J({\Q}_p)$ that generates  $J({\Q}_p)$ as a $p$-adic analytic group. It follows that $X({\Q}_p)$ is not contained in the proper subgroup  $ \Gamma_p$.  Since $X(\Q_p) \cap \Gamma_p$ is a closed proper subvariety of the 1-dimensional compact variety $X(\Q_p)$, it must be compact and  have dimension zero and therefore be finite.  Hence $X({\Q})$ is finite. \vv How can one make this method effective? We describe how this can be done in the next section.
     
     \section{The method of Chabauty as augmented by Coleman}\label{chabauty-coleman}
The Chabauty--Coleman method \cite{Col85a} is one of our most practical tools for actually computing the finite set of rational points on a curve $X$ of genus greater than 1 defined over the rationals,  subject to the same Chabauty condition; namely that the Mordell--Weil rank $r$ of the Jacobian of the curve is strictly less than its genus. \vv

Robert Coleman constructed, in the above conditions, a $p$-adic analytic function $\phi$ on the $p$-adic analytic variety $X^{an}_{/{\Q}_p}$ such that the zeroes of $\phi$ on $X({\Q}_p)$ are\begin{itemize} \item reasonably computable (to any approximation),\item finite in number, and\item include $X({\Q})$.\end{itemize} The construction of such a $\phi$  uses Coleman's $p$-adic abelian integrals on the Jacobian of the curve. \vv

Let $X$ be a curve (of genus $g>1$) defined over the rationals and let $J$ be its Jacobian.  Now fix a prime $p$ {\it of good reduction for} $X$ and a rational point $b \in X(\Q)$. Consider, as before,  the Abel--Jacobi embedding $ i_b: X \rightarrow J$ given by 
$P \mapsto [(P) - (b)].$
Coleman \cite{Col82, Col85b} proved that there is a $p$-adic line integral on holomorphic differentials on the curve satisfying several nice properties (linearity in the integrand, additivity in endpoints, pullbacks under rigid analytic maps, Galois compatibility). The map 
\begin{align*} &J(\Q_p) \times H^0(X_{\Q_p}, \Omega^1) \rightarrow \Q_p \\
&\qquad (Q, \omega)\quad\quad\quad \qquad\mapsto \langle Q, \omega \rangle\end{align*} is additive in $Q$, is $\Q_p$-linear in $\omega$ and is given by
$$\langle Q, \omega \rangle = \langle [D], \omega \rangle \eqqcolon \int_D \omega$$ for $D \in \Div^0(X)$ with $Q = [D]$. Then $$\langle i_b(P), \omega \rangle = \int_b^P \omega.$$
The embedding $i_b$ induces an isomorphism of $g$-dimensional vector spaces $$H^0(J_{\Q_p}, \Omega^1) \simeq H^0(X_{\Q_p}, \Omega^1),$$ giving us the pairing 

\begin{align*}J(\Q_p) \times H^0(J_{\Q_p}, \Omega^1) &\rightarrow \Q_p \\
(Q, \omega_J)\qquad\qquad &\mapsto \int_0^Q \omega_J.\end{align*}
This gives a homomorphism
$$\log: J(\Q_p) \rightarrow H^0(J_{\Q_p}, \Omega^1)^*,$$ {where $\log$ is the logarithm on the $p$-adic Lie group $J(\Q_p)$,} and we have the following diagram

\begin{equation}\label{cc}
\resizebox{8.5cm}{!}{
\begin{tikzpicture}[->,>=stealth',baseline=(current  bounding  box.center)]
 \node[] (X) {$X(\Q)$};
 \node[right of=X, node distance=3.7cm]  (Xp) {$X(\Q_p)$};   
 \node[below of=X, node distance=1.5cm]  (Hf) {$J(\Q)$};
 \node[right of=Hf,node distance=3.7cm] (Hfp) {$J(\Q_p)$};
 \node[right of=Hfp,node distance=3.5cm](Dieu)
  {$H^{0}\left(J_{\Q_p},\Omega^{1}\right)^{\ast} \simeq H^0(X_{\Q_p},\Omega^1)^\ast$}; 

 \path (X)  edge node[left]{\footnotesize } (Hf);
 \path (Xp) edge node[left]{\footnotesize } (Hfp);
 \path (X)  edge (Xp);
 \path (Hf) edge node[above]{} (Hfp);
 \path (Hfp) edge node[above]{$\log $}(Dieu);
 \path (Xp)  edge node[above right]{} (Dieu);
\end{tikzpicture}}
\end{equation}

Recall that under the hypothesis $r < g$, the intersection $$X(\Q_p)_1\coloneqq X(\Q_p) \cap \Gamma_p,$$ and consequently $X(\Q)$, is finite.  Coleman gave a technique to compute $X(\Q_p)_1$, by his construction of $p$-adic integrals that vanish on $\Gamma_p$: in particular, considering an integral of an \emph{annihilating differential} {$\omega$, a holomorphic differential such that $\langle P, \omega \rangle = 0$ for all $P \in J(\Q)$,} then computing the zero locus of this integral on $X(\Q_p)$.  {Bounding the number of zeros of this integral via fairly elementary $p$-adic analysis (for good $p > 2g$) yields the bound $$\#X(\Q) \leq \#X(\Q_p)_1 \leq \#X(\F_p) + 2g-2.$$ In Section \ref{Ex}, we give a worked example of the Chabauty--Coleman method.}

\section{The method of Chabauty--Coleman--Kim}\label{qc}

The construction above crucially uses an assumption that the rank of the Jacobian is small relative to the genus. Nevertheless, there are many interesting curves where this hypothesis is not satisfied, including a number of modular curves we have already seen. \vv
    In a series of papers \cite{kimp1, kim, kimmassey},  Minhyong Kim laid out a program to  extend  Chabauty--Coleman relaxing the condition on Mordell--Weil rank, going beyond the abelian confines of the Jacobian, replacing it by a sequence of \emph{Selmer varieties}, which are carved out of unipotent quotients of $\pi_1^{\et}(X_{\overline{\Q}})_{\Q_p}$, the
   $\Q _p $-\'etale fundamental group of $X_{\overline{\Q}}$ with base point $b$.\vv
   \begin{rem}We will consider {\it unipotent groups} $U$  over $\Q_p$.   These are  algebraic groups over $\Q_p$ possessing a filtration such that the successive graded pieces are isomorphic to the additive group ${\mathbb G}_a$. By a   {\it $\Q_p$-unipotent group}  we will mean  a  topological group isomorphic to $U(\Q_p)$; such a group is a locally compact topological group that admits a filtration where the successive subquotients of the filtration are topological groups isomorphic to finite-dimensional $\Q_p$-vector spaces. By a {\it $\Q_p$-pro-unipotent group} we  mean a projective limit of  $\Q_p$-unipotent groups.
   
    A  {\it $\Q_p$-Malcev completion } ($G\longrightarrow {\bf  G}$) of a topological group $G$ is  the universal solution to the problem of mapping $G$ to  $\Q_p$-pro-unipotent groups.
    
     The {\it $\Q _p $-\'etale fundamental group} is the $\Q_p$-Malcev completion (see  \cite[Appendix A]{hain-matsumoto}) of the usual \'etale fundamental group. \end{rem}
        
    We first recast the Chabauty--Coleman method (see also \cite{Cor19}, \cite{bmAWS}) using $p$-adic Hodge theory, which adds an extra row of compatibilities to diagram \eqref{cc}.  Let $V = \HH^1_{\et}(X_{\overline{\Q}}, \Q_p)^{*}$ and $V_{\dR} \coloneqq \HH^1_{\dR}(X_{\Q_p})^*$, viewed as a filtered vector space with filtration dual to the Hodge filtration.  We have an isomorphism $V_{\dR}/F^0 \simeq \HH^0(X_{\Q_p}, \Omega^1)^*.$ Let $T$ be the set of primes of bad reduction of $X$, together with the prime $p$.  Let $G_{\Q,T}$ be the Galois group of the maximal unramified-outside-$T$ extension of $\Q$.  Let $G_p$ denote the absolute Galois group of $\Q_p$. Then the \'etale formulation of Chabauty--Coleman is given by the following diagram, where the last row is of Bloch--Kato Selmer groups \cite{BK90}:
    
  \begin{equation}\label{eqn:diagram-ab}
\resizebox{8.5cm}{!}{
\begin{tikzpicture}[->,>=stealth',baseline=(current  bounding  box.center)]
 \node[] (X) {$X(\Q)$};
 \node[right of=X, node distance=3.7cm]  (Xp) {$X(\Q_p)$};   
 \node[below of=X, node distance=1.5cm]  (Hf) {$J(\Q)$};
 \node[right of=Hf,node distance=3.7cm] (Hfp) {$J(\Q_p)$};
 \node[right of=Hfp,node distance=3.5cm](Dieu) { $\HH^0(X_{\Q_p},\Omega^1)^\ast$}; 
 
 \node[below of=Hf, node distance=1.5cm] (H1f) {$\HH^1_f(G_{\Q,T},V)$};
 \node[right of=H1f,node distance = 3.7cm] (H1fp){$\HH_f^1(G_{p},V)$};
 \node[right of=H1fp, node distance=3.5cm](H1dR) {$\HH_1^{\dR}(X_{\Q_p})/F^0$};

 \path (X)  edge node[left]{\footnotesize } (Hf);
 \path (Xp) edge node[left]{\footnotesize } (Hfp);
 \path (X)  edge (Xp);
 \path (Hf) edge node[above]{} (Hfp);
 \path (Hfp) edge node[above]{$\log $}(Dieu);
 \path (Xp)  edge node[above right]{$i_b $} (Dieu);
 
 \path (Hf) edge node[above]{} (H1f);
 \path (Hfp) edge node[above]{} (H1fp);
 \path (Dieu) edge node[right] {$\simeq $} (H1dR);
 
\path (H1f) edge node[right]{} (H1fp);
\path (H1fp) edge node[above] {$\simeq $} (H1dR);
\end{tikzpicture}}.
\end{equation}   
  
Let $U_n$ denote the maximal $n$-unipotent
quotient of $\pi_1^{\et}(X_{\overline{\Q}, b})_{\Q_p}$.  We have that $U_1 = V$ and $U_2$ is a central extension

$$1 \xrightarrow{}  \coker(\Q_p(1) \xrightarrow{\cup^*} \wedge^2 V) \xrightarrow{} U_2 \xrightarrow{} V \xrightarrow{} 1.$$

Suppose that $U$ is a Galois-stable quotient of $U_n$.  Kim defined global and local unipotent Kummer maps $j_U$ and $j_{U,v}$ such that the following diagram commutes:
\[
\begin{tikzpicture}[->,>=stealth',baseline=(current  bounding  box.center)]
 \node[] (X) {$X(\Q)$};
 \node[right of=X, node distance=4.7cm]  (Xv) {$\prod_{v \in T}X(\Q_v)$};   
 \node[below of=X, node distance=1.5cm]  (Hf) {$\HH^1(G_{\Q,T},U)$};
 \node[right of=Hf,node distance=4.7cm] (Hfv) {$\prod_{v \in T}\HH^1(G_v,U).$};
\path (X)  edge node[left]{\footnotesize $j_U$} (Hf);
  \path (Xv) edge node[right]{\footnotesize $\prod j_{U,v}$} (Hfv);
 \path (X)  edge (Xv);
 \path (Hf) edge node[above]{\footnotesize $\prod \mathrm{loc}_{v}$} (Hfv);
\end{tikzpicture}
\]

Kim proved that the nonabelian pointed cohomology sets $\HH^1(G_{\Q,T}, U)$ and $\HH^1(G_v, U)$ 
are affine algebraic varieties over $\Q_p$. Motivated by the classical study of Selmer groups, he then refined $\HH^1(G_{\Q,T}, U)$ by local conditions to produce a closed subscheme. Such a closed subscheme is a \emph{Selmer scheme}. We may now give an adapted version \cite{4} of the definition of a \emph{Selmer variety}:

\begin{defn}\label{D:Selmer}
  The  \emph{Selmer variety} $\Sel(U)$ is the reduced scheme associated to the
  subscheme of $\HH^1(G_{\Q,T}, U)$ containing those classes $c$ such that
\begin{itemize}
 \item $\loc_p(c)$ is crystalline,
 \item $\loc_\ell(c) \in j_{U,\ell}(X(\Q_\ell))$ for all $\ell \ne p$,
 \item the projection of $c$ to $\HH^1(G_{\Q,T},V)$ comes from an element of $J(\Q)\otimes
   \Q_p$.
\end{itemize}
\end{defn}\vv

Now the Selmer variety gives rise to the following interesting set of points \[
  X(\Q_p)_{U} \coloneqq  j_{p}^{-1}\left(\loc_{p} (\Sel(U))\right)  \subset X(\Q_p )
\]
and note that  \begin{equation}\label{eq:UUn}
X(\Q) \subset X(\Q_p)_n \coloneqq X(\Q_p)_{U_n}\subset X(\Q_p)_U \,.
\end{equation}\vv

The set $X(\Q_p)_n$ can be computed in terms of $n$-fold iterated Coleman integrals, and one has a series of refinements
$$X(\Q)   \subset \cdots \subset X(\Q_p)_n \subset X(\Q_p)_{n-1} \subset \cdots \subset X(\Q_p)_2 \subset X(\Q_p)_1.$$
Note that the set $X(\Q_p)_1$ is the Chabauty--Coleman set from before. We refer to the points in $X(\Q_p)_n$ as the set of \textbf{Selmer points of level $n$}. We call the points in $X(\Q_p)_n \setminus X(\Q)$ the set of \textbf{mock-rational Selmer points of level $n$}. {Kim has conjectured that for $n \gg 0$, the set $X(\Q_p)_n$ is finite. This conjecture is implied by the conjecture of Bloch--Kato \cite{BK90}.}\vv

Putting everything together, Kim's program is to study finiteness of $X(\Q_p)_U$ using $p$-adic Hodge theory and the following diagram is the nonabelian generalization of \eqref{cc}:
\begin{equation*}
\begin{tikzpicture}[->,>=stealth',baseline=(current  bounding  box.center)]
  \matrix (m) [matrix of math nodes, row sep=3em, column sep=4em, minimum width=2em]
  {
    X(\Q) & X(\Q_p) \\
    \Sel(U) & \HH_f^1(G_{p},U) & U^{\dR}/\textrm{Fil}^0\\
  };
  
  \path[-stealth]
    (m-1-1)   edge (m-1-2)
              edge node [left] {$j_U$} (m-2-1)
    (m-1-2)   edge node [left] {$j_{U,p}$} (m-2-2)
              edge node [above,right] {$\quad j_U^{\dR}$} (m-2-3)
    (m-2-1)   edge node [above] {$\text{loc}_{U,p}$} (m-2-2)
    (m-2-2)   edge node [above] {$\simeq$} (m-2-3);
\end{tikzpicture}
\end{equation*}\vv

Computing the depth-2 Selmer set (or a slightly larger finite set containing it), known as \emph{quadratic Chabauty}, has seen progress in recent years \cite{4, 4.2, 2, 3}, via aspects of the theory of $p$-adic height functions \cite{colemangross, nekovar}. The quadratic Chabauty set $X(\Q_p)_2$ is finite for those curves that satisfy the rank bound \cite{4}
$$r < g + \rho - 1,$$ where {$\rho\coloneqq \rk(\NS(J))$} is the N\'eron--Severi rank of the Jacobian over $\Q$.   {To carry out the quadratic Chabauty method, one uses a nontrivial element of $\ker(\NS(J) \rightarrow \NS(X))$  to construct a nonabelian quotient $U$ of $U_2$, which is used to compute $X(\Q_p)_U$.} \vv

Samir Siksek  \cite{siksek} showed that modular curves of genus 3 or more have $\rho$ at least 2, and consequently, for these curves, quadratic Chabauty allows one to consider Jacobians of higher rank than allowed by Chabauty--Coleman. Balakrishnan, Dogra, M\"uller, Tuitman, and Vonk \cite{4.2, 2} made various aspects of quadratic Chabauty computationally practical, using explicit $p$-adic cohomology to compute {a certain (global) $p$-adic height of {Nekov\'a\v{r}} \cite{nekovar}, depending on a choice of a nontrivial element of $\ker(\NS(J) \rightarrow \NS(X))$.} \vv

Roughly speaking, the method starts from the following observation: the global $p$-adic height admits a decomposition as a sum of local heights: a local height at $p$ that can be computed using $p$-adic Hodge theory, and a finite sum of local heights away from $p$ that, in certain favorable conditions, can be shown to be trivial---or if not trivial, at least  a quantity that can be computed from the geometry of a regular model of the curve. We will discuss this in more detail in the case of bielliptic genus 2 curves in Section \ref{quadptsbiellg2}.\vv

Moreover, the global $p$-adic height is a quadratic form on $H^0(X, \Omega^1)^*$.  Choosing an explicit basis for the space of quadratic forms in terms of Coleman integrals, and knowing sufficiently many rational points (either on $X$ or on $J$) and their $p$-adic heights, one can compute a locally analytic function whose zero locus contains $X(\Q_p)_2$.
\vv

Recently, various new perspectives on the Chabauty--Coleman method and quadratic Chabauty have emerged:
``symmetric Chabauty'' introduced by Siksek  \cite{sikseksymm},  a geometric quadratic Chabauty method introduced by Edixhoven--Lido \cite{EL},  and the $p$-adic Arakelov theoretic one of Besser--M\"uller--Srinivasan \cite{6.1}:\vv
\begin{itemize} \item  Siksek gave a symmetric Chabauty method \cite{sikseksymm}, a variant of the Chabauty--Coleman method (see Section \ref{chabauty-coleman}) for symmetric powers of curves.  Symmetric Chabauty has been used and extended in various ways to determine quadratic points on numerous modular curves $X_0(N)$ \cite{15, 7, NV, 1.1a}.  Box, Gajovi\'c, and Goodman \cite{bgg} further developed a ``partially relative'' symmetric Chabauty method to study cubic and quartic points on several modular curves $X_0(N)$.
 \item   In \emph{geometric quadratic Chabauty}, Edixhoven and Lido \cite{EL} use line
bundles over the Jacobian, the Poincar\'e torsor (a biextension of the Jacobian by ${\mathbb G}_m$), and models over the integers to study rational points under the same rank bound hypothesis.   
 \item Besser, M\"uller, and Srinivasan \cite{6.1} give a new construction of $p$-adic heights on varieties over number fields using $p$-adic adelic metrics on line bundles, in the spirit of Zhang's work on real-valued heights using adelic metrics \cite{zhang}. This leads them to formulate \emph{$p$-adic Arakelov quadratic Chabauty}.
\end{itemize} 
We will not discuss these works in detail in this survey.

\bigskip

\section{Rational points on $X_0(37)$: three perspectives}\label{Ex} 

{As a concrete application of the techniques discussed so far, we present here three perspectives on rational points on the modular curve {\href{https://www.lmfdb.org/Genus2Curve/Q/1369/a/50653/1}{$X_0(37)$}. For further discussion, see Section 5 of \cite{MSw}; and for more, see Section 5 of \cite{7}.}

\vv

  The modular curve $X \coloneqq X_0(37)$ is of genus $2$ and therefore is hyperelliptic. Denote by $$X_0(37)\stackrel{\sigma}{\longrightarrow} X_0(37)$$  its hyperelliptic involution, and 
  by
  $$X_0(37)\stackrel{w}{\longrightarrow} X_0(37)$$  its  Atkin--Lehner involution.  The involutions  $\sigma$ and $w$ commute, generating a Klein group ${\mathcal G}$ of automorphisms. The automorphisms $1, w, \sigma, w\sigma$ are defined over ${\Q}$ and are the only automorphisms of $X_0(37)$ over $\C$.  (See \cite{14.4}.)
\vv
  Form the quotients 
 
  \begin{equation}\label{diagquot}\xymatrix{ \ & X_0(37)\ar[dl]_{i_0}\ar[d]^{x}\ar[dr]^{i_1} & \    \\
  E_0 \coloneqq X_0(37)/\langle \sigma\cdot w\rangle &  \mathbb{P}^1_{/{\Q}} \simeq X_0(37)/\langle \sigma \rangle & E_1\coloneqq  X_0(37)/ \langle w  \rangle  }.\end{equation}

Each $E_j$ is a genus one curve, and we consider each as an elliptic curve with identity element $i_j(\infty)$.  By the Riemann--Hurwitz formula, the ramification locus of each of the double covers:
\begin{equation}\label{d0}X_0(37) \stackrel{i_0}{\longrightarrow}E_0\end{equation}
\centerline{and}
\begin{equation}\label{d1}X_0(37) \stackrel{i_1}{ \longrightarrow} E_1\end{equation}
are ${\Q}$-rational (effective) divisors of degree two:\vv
\begin{itemize}  \item$D_0 \coloneqq \{\eta_0,{\bar{\eta}}_0\} \subset  X_0(37)$---for (\ref{d0}).  \item   $D_1 \coloneqq \{\eta_1,{\bar{\eta}}_1\} \subset  X_0(37)$---for (\ref{d1}).
\end{itemize}
\vv

In particular, $D_1$ is the fixed point set  of $w$  and   $D_0$ is the fixed point set  of $w\sigma$. 

\bigskip

Note that (since $\sigma$ commutes with $w$) each of these involutions  ($\sigma, w, w\sigma$) preserves $D_1$ and $D_0$.
The involution $w\sigma$ interchanges the points  $\eta_1,{\bar{\eta}}_1$. So their image  $e_0\in E_0$---which is therefore the image of a ${\Q}$-rational divisor in $X_0(37)$---is ${\Q}$-rational.  Consequently  $\{\eta_1,{\bar{\eta}}_1\}$ either consists of  a pair of ${\Q}$-rational points{\footnote{That's not the case; see Lemma \ref{fin} below.}} or a conjugate pair of quadratic points in $X_0(37)$.  

For the same reason  the involution $w$ preserves the ramification divisor of $w\sigma$ and interchanges the points  $\eta_0,{\bar{\eta}}_0$ and therefore their image  $e_1\in E_1$ is ${\Q}$-rational. 
\vv
  A visit to the \href{https://www.lmfdb.org/Genus2Curve/Q/1369/a/50653/1}{L-Functions and Modular Forms Database}\footnote{See the sidebar of Related Objects on \href{https://www.lmfdb.org/Genus2Curve/Q/1369/a/50653/1}{LMFDB}.} (LMFDB) \cite{LMFDB} will get you that we have:
\begin{itemize}\label{data}
\item $E_0$ is the elliptic curve $$\href{https://www.lmfdb.org/EllipticCurve/Q/37/b/2}{37.b2}: \ y^2+y=x^3+x^2-23x-50.$$    Its Mordell--Weil group is of order $3$.\vv \item  $E_1$ is the elliptic curve  \begin{equation}\label{37.a1}\href{https://www.lmfdb.org/EllipticCurve/Q/37/a/1}{37.a1}:\ y^2+y=x^3-x.\end{equation}  It has Mordell--Weil rank $1$, and its group of ${\Q}$-rational points is isomorphic to $\Z$. \end{itemize}\vv
  \subsection{Chabauty--Coleman gives finiteness} Let $J_0(37)$ denote the Jacobian of $X_0(37)$.  We have: 

 \[\begin{tikzcd}
	{X_0(37)} & {J_0(37)} \\
	& {E_0 \times E_1}
	\arrow[hook, from=1-1, to=1-2]
	\arrow["\phi", from=1-2, to=2-2]
	\arrow["{i_0 \times i_1}"', from=1-1, to=2-2]
\end{tikzcd}\] 
  
  \vv
  \noindent where $i_0,i_1$ are (as above) the modular parametrization of $E_0,E_1$, and $$\phi: J_0(37)  \onto E_0\times E_1$$ is an isogeny. 
  Since $ \{E_0\times E_1\}({\Q})$  is---by the data above---a group isomorphic to ${\Z}\times{\Z}/3{\Z}$ (contained in  a cyclic group of order three   times the elliptic curve  $E_1$) we see that  the Zariski closure of the group of ${\Q}$-rational points $J_0(37)({\Q})$ is an algebraic subgroup in  $J_0(37)$ of codimension $1$  so can intersect only finitely with $X_0(37)$---giving that  $X_0(37)(\Q)$ is finite. \vv
  
  \subsection{The projection to $E_0$ gets us the precise set of ${\Q}$-rational points}
  The cusp ${\bf \infty}\in  X_0(37)$ is a  ${\Q}$-rational point, as are the four points  
\begin{equation}\label{pts}{\mathcal S}={\mathcal G}\cdot{\bf \infty}=\{ {\bf \infty},\  w({\bf \infty})={\rm the\ cusp\ }{\bf 0}, \ \sigma({\bf \infty}),  \sigma({\bf 0}) \}.\end{equation}
  \vv
\begin{thm}  These are the only four ${\Q}$-rational points on $X_0(37)$. \end{thm}\vv

  \begin{proof}   Returning to the mapping of degree two 
      $$X_0(37)({\Q}) \stackrel{i_0}{\longrightarrow} E_0({\Q}),$$  
     
     \noindent since $E_0({\Q})$ is cyclic of order three, we see that
     
     \begin{itemize}\item the pair $\{ {\bf \infty}, \sigma w({\bf \infty})=\sigma({\bf 0})\}$  maps to the origin in $ E_0({\Q})$  and 
     \item the pair $\{{\bf 0},   \sigma w({\bf 0})= \sigma({\bf \infty})\}$ maps to a (nonzero) point $e\in   E_0({\Q})$. 
     \item  Recalling that the pair $\{\eta_1,{\bar{\eta}}_1\}$ discussed above maps to $e_0\in E_0({\Q})$ and noting that $e_0$ cannot be any of the above two ${\Q}$-rational points of $E_0$, it must be the third  ${\Q}$-rational point, giving us:
     $$ e_0 = 2e=-e \in E_0({\Q}).$$\end{itemize}\vv
       These three bullet points show that $$i_0^{-1}(E_0(\Q)) =  {\mathcal{S}} \cup \{\eta_1,{\bar{\eta}}_1\}.$$  Since $$X_0(37)({\Q})  \subset i_0^{-1}(E_0({\Q}))$$\vv 
       \noindent it is enough to show that $\eta_1,{\bar{\eta}}_1  \notin X_0(37)({\Q})$.  This is established in Lemma \ref{fin} below, completing the proofs.\end{proof} \vv

\begin{lem}\label{fin} The pair of points  $\{\eta_1, {\bar \eta}_1\} \in X_0(37)({\bar{\Q}})$   are  ${\Q}({\sqrt{37}})$-conjugate points   defined over $\Q({\sqrt{37}})$ and not over $\Q$. \end{lem} \vv

   {\bf Proof of Lemma \ref{fin}} \vv
   The involutions $\sigma$ and $w$ of $X_0(37)$ are easily described in terms of  the model \footnote{Note that the LMFDB gives a different ``simplified'' model of $X_0(37)$ which can be seen to be isomorphic to ours, using the following  \texttt{Magma} commands:\\
\texttt{> R<x> := PolynomialRing(RationalField());}\\
\texttt{> X37v1 := HyperellipticCurve(x\^\,6 + 8*x\^\,5 - 20*x\^\,4 + 28*x\^\,3 - 24*x\^\,2 +12*x - 4);} \\
\texttt{> X37v2 := HyperellipticCurve(-x\^\,6 - 9*x\^\,4 - 11*x\^\,2 + 37);}\\
\texttt{> IsIsomorphic(X37v1,X37v2);}\\
\texttt{true Mapping from: CrvHyp: X37v1 to CrvHyp: X37v2}\\
\texttt{with equations : }\\
\texttt{1/2*\$.1 - \$.3}\\
\texttt{1/2*\$.2}\\
\texttt{1/2*\$.1}\\
\texttt{and inverse}\\
\texttt{\$.3}\\
\texttt{1/4*\$.2}\\
\texttt{-1/2*\$.1 + 1/2*\$.3}}
of $X_0(37)$ given by 
   \begin{equation}\label{model}y^2 = g(x) \coloneqq -x^6 -9x^4 - 11x^2 + 37,\end{equation} as found by Mazur and Swinnerton-Dyer \cite[\S 5.1]{MSw}. We have: \begin{enumerate}\label{invs} \item[(a)] $(x,y)\ \stackrel{\sigma}{\mapsto}\ (x,-y)$, 
   \item[(b)] $(x,y)\ \stackrel{w}{\mapsto}\ (-x,y)$  and 
 \item[(c)] $(x,y)\ \stackrel{w\sigma }{\mapsto}\  (-x,-y).$\end{enumerate} 

\vv
   The proof  of (c)     follows  from (a) and (b) by composition.
\vv
   The proof  of (a) is  simply that the quotient  of the involution $(x,y)\  {\mapsto}\ (x,-y)$ is  of genus zero as is clear from the equation; so that involution is the hyperelliptic involution $\sigma$.
   
 \vv
 The proof of (b) follows from considering  the following model\footnote{This alternative model of $E_1$ can be checked to be isomorphic to \eqref{37.a1}, using the following \texttt{Magma} commands: \\
\texttt{> P<u,v,z> := ProjectiveSpace(Rationals(),2);}\\
\texttt{> C := Curve(P, v\^\,2*z - (-u\^\,3-9*u\^\,2*z-11*u*z\^\,2+37*z\^\,3));}\\
\texttt{> pt := C![1,4,1];}\\
\texttt{> E := EllipticCurve(C,pt);}\\
\texttt{> E1 := EllipticCurve([0,0,1,-1,0]);}\\
\texttt{> IsIsomorphic(E1,E);}\\
\texttt{true}} (over ${\Q}$) for  the expression  of  $X_0(37)$ as the double cover $X_0(37) \stackrel{i_1}{\longrightarrow} E_1 = X_0(37)^+$ over ${\Q}$:

 \begin{equation}\label{dc}\xymatrix{X_0(37):\ar[d]^{i_1}  &  y^2 =  -x^6 -9x^4 - 11x^2 + 37\\
 E_1:&  v^2 =  -u^3 -9u^2 - 11u + 37\ar[u]^{u=x^2;\   v=y}}\end{equation}
 \vv
 Since  $ \{\eta_1,{\bar{\eta}}_1\}$ consists of the fixed points of  the involution $w$, we have:
 $$ \{\eta_1,{\bar{\eta}}_1\} = \{ (0, \pm {\sqrt{37}})\},$$

 \noindent from which it follows that $$w\sigma:(0, \pm {\sqrt{37}})\ \mapsto\ (0, \mp {\sqrt{37}})$$
 and therefore $$i_0(\eta_1)=i_0({\bar{\eta}}_1) = i_0(0, +{\sqrt{37}}) = i_0(0, -{\sqrt{37}})\ \in \ E_0(\Q),$$  i.e., it is a ${\Q}$-rational point of $E_0$, which can be neither $i_0({\bf \infty})$ nor $i_0({\bf 0})$ so must be the third $\Q$-rational point.\vv  
  
 \subsection{The Chabauty--Coleman method would also give us the set of ${\Q}$-rational points}
Since $J_0(37) \sim E_0 \times E_1$, we have $$\rk\; J_0(37)(\Q) = \rk\; E_0(\Q) + \rk\;E_1(\Q) = 1,$$ where the ranks of the elliptic curves $E_0, E_1$ were earlier found in LMFDB. Since the rank of $J_0(37)$ is less than the genus of the curve $X\coloneqq X_0(37)$, we we may carry out the Chabauty--Coleman method on \eqref{model} to compute the set $X(\Q)$. We use the prime $p = 3$ and take $\{\frac{dx}{y}, \frac{x \,dx}{y}\}$ as our basis of $H^0(X,\Omega^1).$

\vv

Searching in a box for rational points of small height, one finds the points $(\pm 1, \pm 4) \in X(\Q)$.  The point $P \coloneqq [(1,-4) - (-1,4)] \in J_0(37)(\Q)$ is non-torsion, since the 3-adic Coleman integral  {of a holomorphic differential along this point} is nonzero: $$\int_{(-1,4)}^{(1,-4)} \frac{x\, dx}{y} = 3^2 + 2 \cdot 3^3 + 3^4 + 2 \cdot 3^5 + 3^7 + O(3^9).$$ Moreover, $\int_{(-1,4)}^{(1,-4)} \frac{dx}{y} = O(3^9)$.  In fact, $\int_{(-1,4)}^{(1,-4)} \frac{dx}{y}$ is identically 0, as can be seen by applying the involution $w\sigma$ to the integrand and the endpoints of involution.  Thus we may take $\frac{dx}{y}$ as our annihilating differential. The curve $X$ over $\F_3$ has the following rational points:
$$(\overline{0},\overline{1}), (\overline{0},\overline{2}), (\overline{1},\overline{1}), (\overline{1},\overline{2}), (\overline{2},\overline{1}), (\overline{2},\overline{2}) \in X_{\F_3}(\F_3),$$
which correspond to the residue disks over which we carry out our computation. 

\vv

Fixing as our basepoint $(-1,4) \in X(\Q)$, we start in the residue disk corresponding to $(\overline{0},\overline{1})$. We take the point  $S_0 = (0, y_0)$ in the residue disk, where $y_0 \in \Z_3$ is the unique square root of 37 that satisfies $y_0 \equiv 1 \pmod{3}$. We compute our local coordinate at $S_0$: since $x = 0$, we take $x(t) = t$. Then $y(t)$ correspondingly is found by applying Hensel's lemma to produce the $3$-adic power series expansion of the square root \begin{align*}y(t) &=  \sqrt{g(x(t))},\\
&= \sqrt{-t^6 -9t^4 - 11t^2 + 37},\end{align*}
which yields the following:
\begin{align*}S_t 
&= (t,     -3788  + (2159 + O(3^{10}))t^2   - (15737 +  O(3^{10}))t^4 + \\
&\qquad - (23833 + O(3^{10}))t^6  +  (746\cdot 3^3 + O(3^{10}))t^8 + O(t^{10})) \\
&=: (x(t),y(t)). \end{align*} \vv

         We now wish to compute the zeros of the power series $I(3T)$, where
$$I(T) = \int_{(-1,4)}^{S_0} \frac{dx}{y} + \int_{S_0}^{S_T} \frac{dx(t) \, dt}{y(t)}.$$
Using the 3-adic power series expansions calculated for $x(t),y(t)$ above, as well as the value of the Coleman integral between $(-1,4)$ and $S_0$, this yields 
\begin{align*}I(3T) &=  \left(3 + 3^{3} + 2 \cdot 3^{4} + 2 \cdot 3^{5} + 3^{6} + 3^{7} + 2 \cdot 3^{8} + 3^{9} + 3^{10} + O(3^{11})\right)T +\\
& \left(3^{2} + 2 \cdot 3^{4} + 2 \cdot 3^{5} + 3^{7} + 2 \cdot 3^{8} + 2 \cdot 3^{9} + 3^{10} + O(3^{12})\right)T^{3} + \\
& \left(3^{6} + 3^{7} + 2 \cdot 3^{8} + 3^{9} + 3^{10} + 3^{11} + 2 \cdot 3^{13} + 2 \cdot 3^{14} + O(3^{15})\right)T^{5} +\\
&  \left(3^{8} + 2 \cdot 3^{9} + 3^{10} + 2 \cdot 3^{11} + 2 \cdot 3^{12} + 2 \cdot 3^{13} + 2 \cdot 3^{15} + O(3^{17})\right)T^{7} + \\
& \left(3^{7} + 2 \cdot 3^{8} + 2 \cdot 3^{10} + 2 \cdot 3^{11} + 3^{12} + 3^{14} + 2 \cdot 3^{16} + O(3^{17})\right)T^{9} + O(T^{10}),\end{align*}
which has precisely one zero at $T = 0$, corresponding to $S_0$, which we can identify as $(0, \sqrt{37})$.  \vv

Continuing in this way, parametrizing each residue disk by a local coordinate and computing the zeros of the corresponding $I(3T)$ in each residue disk, we find that
$$X(\Q_3)_1 = \{(0, \pm \sqrt{37}), (\pm 1, \pm 4)\},$$ from which we immediately produce $$X(\Q) = (\pm 1, \pm 4).$$ \vv

\begin{rem}It was fairly lucky that $X(\Q_3)_1 = \{(0, \pm \sqrt{37}), (\pm 1, \pm 4)\}$ and was not much larger.  {Finding a small good prime $p$ such that there are no  mock-rational Selmer points---or where the mock-rational points are easily-recognized algebraic points---may be an issue.}  By the Weil bound, we know that $\#X(\F_p)$ grows linearly as $p$ grows.  So if we had used a larger prime $p$ in the Chabauty--Coleman method, we would expect more $p$-adic points in $X(\Q_p)_1$, and we may not be able to immediately recognize these extra points.\end{rem}
  
\section{Quadratic points on bielliptic curves of genus $2$ using quadratic Chabauty}\label{quadptsbiellg2}
{In the previous section, we considered the problem of determining the finitely many rational points on $X_0(37)$. We could also study the finite sets $X_0(37)(K)$ for various other number fields $K$, one number field at a time.  Or, we could further study $\Symm^d(X_0(37))(\Q)$, as described in Part \ref{overview}, which would tell us about all degree $d$ points on $X_0(37)$.}
 
\vv

{We start by considering $X_0(37)(K)$ for a fixed quadratic field $K$. If the rank of $J_0(37)(K)$ is now 2, and if this is because the rank of $E_0(K)$ increases to 1---recall from Section \ref{Ex} that the rank of $E_0(\Q)$ is 0---then the Chabauty--Coleman method no longer applies. However, since $X_0(37)$ is bielliptic and genus 2, we can use the method of \cite{4}, which gives a particularly explicit description of quadratic Chabauty functions (an elaboration of the case $n = 2$ in Section \ref{qc})  using $p$-adic height functions and double Coleman integrals on elliptic curves, for bielliptic genus 2 curves.}   We describe this below in some generality, and then use it to study rational points on $X_0(37)$ over $K=\Q(i)$.

\vv

Let $K = \Q$ or a quadratic imaginary extension, and let $X/K$ be a genus 2 bielliptic curve $$y^2 = x^6 + a_4x^4 + a_2x^2 + a_0,$$ with $a_i \in K$. Let $C_1$ and $C_2$ be the elliptic curves over $K$ defined by the equations
$$C_1: y^2 = x^3 + a_4 x^2 + a_2 x + a_0, \qquad C_2: y^2 = x^3 + a_2x^2 + a_4a_0 x + a_0^2,$$
and let $f_1: X \rightarrow C_1$ be the map that sends $(x,y)$ to $(x^2,y)$ and $f_2: X \rightarrow C_2$ be the map that sends $(x,y) \rightarrow (a_0 x^{-2}, a_0 yx^{-3})$. 
\vv

 We will be considering the case where the Mordell--Weil ranks of $C_1$ and $C_2$ over $K$ are equal to $1$. Letting $J$ denote the Jacobian of $X$ we have that the rank of $J$  over $K$ is $2$. The natural mapping defined over $K$   \begin{equation}\label{Sym2}  \Symm^2(X) \to J \end{equation} (i.e., setting $p=2$ in Equation \ref{Sym})
  is\begin{itemize} \item  an isomorphism if $X$ is not hyperelliptic, or is\item  an isomorphism in the complement of an exceptional fiber ${\mathcal E}\subset \Symm^2(X)$ isomorphic to ${\mathbb P}^1$ over $K$ if $X$ is hyperelliptic.\end{itemize}
 In other words, all quadratic points of $X$ over $K$ are parameterized, in an appropriate sense, either by $J(K)$, if $X$ is not hyperelliptic, or by $J(K)$ together with the isomorphism
 ${\mathcal E}\simeq {\mathbb P}^1(K)$, if $X$ is hyperelliptic. \vv

   These parametrizations are neat, and explicit, but they still leave untouched the question: for a given quadratic field $K$ what---exactly---is the finite set $X(K)$?  We want to use quadratic Chabauty to answer such questions.
   
   \vv {Fix some auxiliary choices, including an id\`ele class character $\chi: G_K^{\textrm{ab}} \rightarrow \Q_p$.} Define $Z_1, Z_2 \subset X \times X$ to be the graphs of the automorphisms $g_1: (x,y) \mapsto (-x,y)$ and $g_2: (x,y) \mapsto (-x,-y)$, respectively. We take as our correspondence $Z \coloneqq Z_1 - Z_2$.  When $K = \Q$ fix a prime $\fp = (p)$ to be a prime of good ordinary reduction. When $K$ is imaginary quadratic, take $p$ to be a rational prime that splits as $\fp\overline{\fp}$ where both $\fp$ and $\overline{\fp}$ are primes of good ordinary reduction. \vv

Let $h_{C_1}$ and $h_{C_2}$ denote the global $\fp$-adic height functions associated to the choices made above and $h_{C_i, \fp}$ the respective local height at $\fp$, with the global height written as the sum of local heights $$h_{C_i} = \sum_v h_{C_i, v}.$$ Suppose $C_1(K)$ and $C_2(K)$ each have Mordell--Weil rank 1, and let $P_i \in C_i(K)$ be points of infinite order.
Let $$\alpha_i = \frac{h_{C_i}(P_i)}{[K:\Q]\log_{C_i}(P_i)^2}.$$  Let $\Omega$ denote the finite set of values taken by $$-\sum_{v\nmid p} (h_{C_1, v}(f_1(z_v)) - h_{C_2,v}(f_2(z_v)) - 2\chi_v(x(z_v))),$$ for $(z_v) \in \prod_{v\nmid p} X(K_v).$ Then $X(K)$ is contained in the finite set of $z \in X(K_{\fp})$ cut out by the quadratic Chabauty function $$h_{C_1,\fp}(f_1(z)) - h_{C_2, \fp}(f_2(z)) - 2\chi_{\fp}(x(z)) - \alpha_1\log_{C_1}(f_1(z))^2 + \alpha_2\log_{C_2}(f_2(z))^2 \in \Omega,$$ where  $$\log_{C_i}(Q) = \int_{\infty}^Q \frac{dx}{2y},$$ the single Coleman integral we saw in the Chabauty--Coleman method (with $\infty$ denoting the point at infinity on the corresponding elliptic curve) and $$h_{C_i,\fp}(z)$$ is a double Coleman integral.

\vv
\begin{rem}Over $K = \Q(i)$, the elliptic curves \href{https://www.lmfdb.org/EllipticCurve/Q/37/a/1}{37.a1} and \href{https://www.lmfdb.org/EllipticCurve/Q/37/b/2}{37.b2} each  have rank 1.  The computation  in \cite{4}, applies quadratic Chabauty as described above at the primes $p = 41, 73, 101$ {to produce, for each prime $p$, a finite superset of $p$-adic points containing $X(K)$. This is then combined with another method, the \emph{Mordell--Weil sieve}}, to give \begin{equation}\label{Qi}X_0(37)(K)= \{(\pm 2i, \pm 1),(\pm 1,  \pm 4), {\bf \infty},  {\bf 0}\}.\end{equation} \end{rem}
\vv
\subsection{Explicitly determining quadratic points} 
Quadratic Chabauty for bielliptic curves over $\Q$ was subsequently refined by Bianchi \cite{bianchi}  using $p$-adic sigma functions in place of double Coleman integrals. This was recently extended by Bianchi and Padurariu \cite{bianchi-padurariu}, where an implementation was given to study rational points on {\it all} rank 2 genus 2 bielliptic curves in the LMFDB, including the Atkin--Lehner quotient curve  $X_0(166)^* \coloneqq X_0(166)/\langle w_2, w_{83} \rangle$ (with LMFDB label \href{https://www.lmfdb.org/Genus2Curve/Q/13778/a/27556/1}{13778.a.27556.1}),  as well as the Shimura curve $X_0(10,19)/\langle w_{190} \rangle$.

\vv
Using a slight extension of their work to $K = \Q(i)$, as done in \cite{BMcode}, one can use a smaller prime to carry out the computation of a finite set containing the depth 2 Selmer set for $X_0(37)$. (Recall  Definition \ref{D:Selmer} in   Section \ref{qc}.)    {We carried out this computation for $p=13$ and recovered the points $(\pm 2i, \pm 1),(\pm 1,  \pm 4),$ and ${\bf \infty},  {\bf 0}$.  But lurking within the set of depth 2 Selmer points, we also found the \emph{algebraic}  points $(\pm \sqrt{-3}, \pm 4),$ these being initially observed 73-adically in \cite{4}.  We also found several other  mock-rational Selmer points, such as 
\small{$$(5 + 8 \cdot 13 + 12 \cdot 13^{2} + 4 \cdot 13^{3} + 2 \cdot 13^{4} + 3 \cdot 13^{5}  +  O(13^{6}), 1 + 3 \cdot 13 + 3 \cdot 13^{2} + 9 \cdot 13^{3} + 12 \cdot 13^{4} + 5 \cdot 13^{5}  + O(13^{6})).$$}}  
  
 {
See   Banwait--Najman--Padurariu \cite{BNP} for an extensive discussion---and for results---regarding  quadratic points on $X_0(N)$. In particular they show that {$X_0(37)(\mathbb{Q}(\sqrt{d})) = X_0(37)(\Q)$} for $$d= -6846, -2289, 213, 834, 1545, 1885, 1923,
2517, 2847, 4569,$$
$$ 6537, 7131, 7302,
7319, 7635, 7890, 8383, 9563, 9903.$$
  
  {We could continue by varying the quadratic fields $K$, and in principle, if the rank is not too large, apply Chabauty--Coleman, quadratic Chabauty or variations thereof---possibly combining with other Diophantine techniques---to determine the $K$-rational points on $X_0(37)$. But eventually the ranks outpace our current collection of Diophantine tools. For instance, over $K = \Q(\sqrt{-139})$, a \texttt{Magma} computation reveals that the elliptic curve $E_0$ has rank 3, as does $E_1$, and so $J_0(37)(K)$ here altogether has rank 6, making it a challenge for existing methods.  
  
  \vv
  
Now indeed, since $X_0(37)$ is hyperelliptic, it has \emph{infinitely many} quadratic points. Nevertheless, one can describe all quadratic points on $X_0(37)$, using the $\Symm^2$ perspective and the maps to the various quotients of $X_0(37)$ in the diagram \eqref{diagquot}, as was done by Box \cite{7}. The hyperelliptic covering map $x:  X_0(37) \rightarrow \mathbb{P}^1$ is one source of infinitely many rational points, and the rank 1 elliptic curve quotient $E_1$ is another source of infinitely many rational points. Finally, the elliptic curve quotient $E_0$ gives three rational points, and Box pieced together these three sources of rational points to describe $\Symm^2(X_0(37))(\Q)$, as below. \vv

\vv The $x$-map gives us all points $\{(x_i,\sqrt{g(x_i)}), (x_i,-\sqrt{g(x_i)})\} \in \Symm^2(X_0(37))(\Q)$, where $x_i$ ranges through all rational numbers. We can find $P_1 \in X_0(37)(\Q(\sqrt{-3}))$ such that $[P_1 + \overline{P_1} - {\bf \infty}-{\bf  0}]$ generates the free part of the Mordell--Weil group of $J_0(37)(\Q)$, and we have the points $\mathcal{P}_{1,0} \coloneqq \{P_1, \overline{P_1}\}$ and $\mathcal{P}_{0,1} \coloneqq \{{\bf \infty}, w({\bf  0})\}$. Finally, for any $(a,b) \in \Z \times \Z/3\Z \setminus\{(0,0)\}$, there is a point $\mathcal{P}_{a,b} \in \Symm^2(X_0(37)(\Q))$ defined by the unique effective degree 2 divisor $P$ such that $$P - {\bf \infty} - {\bf 0} \sim a\mathcal{P}_{1,0} + b\mathcal{P}_{0,1}  - (a+b)({\bf \infty} + {\bf 0})$$ for any lift of $b$ to $\Z$. }

\section{Thanks}
\vv
 We are grateful to Barinder Banwait, Alexander Betts, Francesca Bianchi, Maarten Derickx, Netan Dogra,  Minhyong Kim, Steffen M\"uller, Filip Najman, Ken Ribet, Matthias Sch\"utt, Bianca Viray, Isabel Vogt, Preston Wake, and the anonymous referees for their illuminating and extremely helpful comments.  Thanks also to Netan Dogra for providing the appendix on Bring's curve. Thanks to the organizers in the IAS for organizing and hosting the Frank C. and Florence S. Ogg Professorship  in Mathematics, and thanks to Andrew for inspiring all of us. The research for this paper was partially supported by NSF grant DMS-1945452, the Clare Boothe Luce Professorship (Henry Luce Foundation), Simons Foundation grant no. 550023 (J.S.B.), and NSF grant DMS-2152149 (B.M.).

\appendix
\section{Quadratic points on Bring's curve, by Netan Dogra}

We consider \emph{Bring's curve}, the smooth projective genus 4 curve $X$ in $\mathbb{P}^4$ given as the locus of common zeros of the following system of equations:
\begin{equation}\label{brings}
\begin{cases}
x_1 + x_2 + x_3 + x_4 + x_5 &= 0 \\
x_1^2 + x_2^2+ x_3^2 +  x_4^2 + x_5^2 &= 0 \\
x_1^3 + x_2^3+ x_3^3 +  x_4^3 + x_5^3  &= 0.  \end{cases}
\end{equation} \vv

From the quadratic defining equation of Bring's curve, we see that $X(\R) = \emptyset$, so we have that $X(\Q) = \emptyset$. However, considering the curve  instead over $K = \Q(i)$, we see several $K$-rational points: for instance, all permutations of the coordinates of the points $(1: \pm i: -1: \mp i: 0)$ are in $X(\Q(i))$. Could there possibly be more points?

\begin{prop} The only quadratic points on Bring's curve are over $\Q(i)$, and up to permutation of coordinates, they are $(1: \pm i: -1: \mp i: 0).$
\end{prop}

The automorphism group of $X$ is the symmetric group $S_5$, given by permutation of the five coordinates.   Using the action of $S_5$ on $X$, one can see that the Jacobian $J$ of $X$ is isogenous to $E^4$ \cite[Section 8.3.2]{serre-galois}, where $E$ is the rank zero elliptic curve with LMFDB label \href{https://www.lmfdb.org/EllipticCurve/Q/50/a/3}{50.a3}. Since Bring's curve is not hyperelliptic, the map $$\Symm^2(X) \hookrightarrow \Pic^0(X)$$ is injective, and  since $\Pic ^0 (X)({\mathbb Q})$ is finite it follows  that there are only finitely many quadratic points on Bring's curve.  \vv

There is also a simple description of a map $\Symm ^2 (X)\to E^4 $ with finite fibers. The quotient of Bring's curve by the involution swapping two coordinates is isomorphic to the curve
\[
E': x^3 +y^3 +1 + x^2 y +y^2 x + x^2 +y^2 +xy+x+y=0
\]
by projecting the three non-permuted coordinates to $\mathbb{P}^2 $.
This is isomorphic to the elliptic curve
\[
E:y^2 +5x^3+5x^2+4 = 0 \qquad\textrm{(LMFDB label \href{https://www.lmfdb.org/EllipticCurve/Q/50/a/3}{50.a3})} 
\] 
via 
\[
(x,y)\mapsto \left( \frac{2}{1+2x+2y},\frac{4(y-x)}{1+2x+2y} \right) .
\]
We have $E(\Q) = \{\infty, (-2,\pm 4)\}$. The $S_3$-action on $E'$ corresponds to the action of $E(\Q )$ and $-1$ on $E$.  \vv

Now fix a quadratic point $P=(x_0 :x_1 :x_2 :x_3 :1 )$ on Bring's curve. Up to an $S_3 $ permutation, we may assume it maps to $\infty$ in $E$ after quotienting by the involution switching $x_0 $ and $x_1 $. Suppose $\sigma$ generates the Galois group of the field of definition of $P$. Let $x=x_2 $ and $y=x_3 $. Then
\[
\frac{y-x}{1+2x+2y}=-\frac{\sigma y-\sigma x}{1+2\sigma x+2\sigma y}.
\]
This reduces to the equation
\[
\Tr y+4\Nm y =\Tr x +4\Nm x. 
\]
Thus quadratic points on Bring's curve are 5-tuples $(x_1 :x_2 :x_3 :x_4 :x_5 )$ of quadratic points in $\mathbb{P}^5$ satisfying, for all ${i_1 ,i_2 ,i_3 }\subset \{1 ,2,3,4,5\} $,
\[
\prod _{\sigma \in S_3 }(\Tr (x_{i_{\sigma (1)}}/x_{i_{\sigma (3)}})+4\Nm (x_{i_{\sigma (1)}}/x_{i_{\sigma (3)}})-\Tr (x_{i_{\sigma (2)}}/x_{i_{\sigma (3)}})-4\Nm (x_{i_{\sigma (2)}}/x_{i_{\sigma (3)}}))=0.
\] \vv

Up to the $S_5$-action, we may reduce to finding tuples $(x_1 ,x_2 ,x_3 ,x_4 )$ defining a quadratic point $(x_1 :x_2 :x_3 :x_4 :1)$ on Bring's curve and satisfying 
\[
\Tr x_1+4\Nm x_1 =\Tr x_2 +4\Nm x_2 . 
\]
and either 
\[
\Tr x_1+4\Nm x_1 =\Tr x_3 +4\Nm x_3 ,
\]
\[
\Tr \left(\frac{1}{x_1}\right) +4\Nm\left(\frac{1}{x_1}\right) =\Tr \left(\frac{x_3}{x_1}\right) +4\Nm \left(\frac{x_3}{x_1}\right),
\]
or
\[
\Tr \left(\frac{1}{x_3}\right) +4\Nm \left(\frac{1}{x_3}\right) =\Tr \left(\frac{x_1}{x_3}\right) +4\Nm \left(\frac{x_1}{x_3}\right).
\]
\vv
Write each quadratic point $x_i = u_i + w_i$, where $u_i$ and $w_i$ are in plus and minus eigenspaces for the Galois involution. These equations define a finite scheme over $\Q $, and one may check that its rational points correspond exactly to the quadratic points in the statement of the proposition.

\ed
\begin{thebibliography}{bib}
\bibitem{Zan} U. Zannier, Torsion in algebraic groups and problems which arise, To appear.
\end{thebibliography}
